\documentclass{article}
\usepackage{amsmath}
\usepackage{amssymb}
\usepackage{ulem}

\evensidemargin=\oddsidemargin
\newdimen\dummy
\dummy=\oddsidemargin
\newtheorem{theorem}{Theorem}[section]

\newtheorem{corollary}[theorem]{Corollaire}

\newtheorem{definition}[theorem]{Définition}

\newtheorem{lemma}[theorem]{Lemme}

\newtheorem{proposition}[theorem]{Proposition}
\newtheorem{remark}[theorem]{Remarque}

\newenvironment{proof}[1][Preuve]{\textbf{#1 :} }{\ \rule{0.5em}{0.5em}}

\begin{document}

\author{Justin FEUTO - Ibrahim FOFANA - Konin KOUA \and \ \ \ \ \ Universit%
\'{e} de Cocody - Laboratoire de Math\'{e}matiques Fondamentales \and 22
B.P. 582 Abidjan 22 \and justfeuto@yahoo.fr, fofana\_ib\_math\_ab@yahoo.fr,\
kouakonin@hotmail.com}
\title{Espaces de fonctions \`{a} moyenne fractionnaire int\'{e}grable sur
les groupes localement compacts.\ }
\date{}
\maketitle

\begin{abstract}
Let $G$\ be a locally compact group which is $\sigma $-compact, endowed with
a left Haar measure $\lambda .$\ Denote by $e$\ the unit element of $G$, and
by $B$ an open relatively compact and symmetric neighbourhood of $e$. For
every $\left( p,q\right) $ belonging to $\left[ 1\ ;\ +\infty \right] ^{2}$,
we give an equivalent and a priori more manageable definition of the Banach
space $L_{(q,p)}^{\pi }(G),$\ defined by R. C. Busby and H. A. Smith in \cite%
{1}. In the case $G$ is a group of homogeneous type, we look at the
subspaces $\left( L^{q},L^{p}\right) ^{\alpha }(G)$\ of the space $%
L_{(q,p)}^{\pi }(G)$. Theses subspaces are extensions to non abelian groups
of the spaces of functions with integrable mean, defined by I. Fofana in 
\cite{2}. Finally we show that $L^{\alpha ,+\infty }(G)$\ \ is a complex
subspace of $\left( L^{q},L^{p}\right) ^{\alpha }(G)$.\newline
\begin{equation*}
\mathbf{Resum\acute{e}} 
\end{equation*}
\newline
Soit $G$ un groupe localement compact \ $\sigma $-compact, muni d'une mesure
de Haar \`{a} gauche $\lambda .$ D\'{e}signons par $e$ l'\'{e}l\'{e}ment
neutre de $G$, et par $B$ un voisinage ouvert relativement compact et sym%
\'{e}trique de $e$. Pour tout couple \ $\left( p,q\right) $\ d'\'{e}l\'{e}%
ments de $\left[ 1\ ;\ +\infty \right] $, nous donnons une d\'{e}finition 
\'{e}quivalente et \`{a} priori plus maniable de l'espace $L_{(q,p)}^{\pi
}(G),$\ d\'{e}fini par R. C. Busby et H. A. Smith dans \cite{1}. Dans le cas
o\`{u} $G$ est un groupe de type homog\`{e}ne, nous \'{e}tudions les
sous-espaces $\left( L^{q},L^{p}\right) ^{\alpha }(G)$\ de l'espace $%
L_{(q,p)}^{\pi }(G)$. Ces sous-espaces sont des extensions aux groupes non ab%
\'{e}liens, des espaces de fonctions \`{a} moyenne fractionnaire int\'{e}%
grable, d\'{e}finis par I. Fofana dans \cite{2}. Nous montrons par la m\^{e}%
me occasion que $L^{\alpha ,+\infty }(G)$\ est un sous-espace vectoriel
complexe de $\left( L^{q},L^{p}\right) ^{\alpha }(G)$.
\end{abstract}

\newpage

\section{ NOTATIONS ET RESULTATS}

$G$\ est un groupe localement compact, $\sigma -$compact, d'\'{e}l\'{e}ment
neutre $e,$ muni d'une mesure de Haar \`{a} gauche not\'{e}e $\lambda .$

$L_{0}(G)$\ d\'{e}signe l'espace vectoriel complexe des classes d'\'{e}%
quivalences modulo l'\'{e}galit\'{e} $\lambda -$presque partout des
fonctions complexes $\lambda $-mesurables sur $G$. Pour tout sous-ensemble $%
E $\ de $G$, $\chi _{_{E}}$\ d\'{e}signe la fonction caract\'{e}ristique de $%
E$.

Pour tout \'{e}l\'{e}ment $p$\ de $\left[ 1\ ;\ +\infty \right] $, $\left\|
{}\right\| _{p}$\ d\'{e}signe la norme usuelle dans l'espace de Lebesgue $%
L^{p}(G).$

a) Consid\'{e}rons une $U-V$\ partition uniforme $\pi $\ de $G$\ ( voir D%
\'{e}finition \ref{DEF: Part-uni}) et deux \'{e}l\'{e}ments $p$\ et $q$\ de $%
\left[ 1\ ;\ +\infty \right] .$

Posons :

$\qquad \blacktriangleright \ $Pour tout \'{e}l\'{e}ment $f$\ de $L_{0}(G),$

\begin{equation*}
\left\| f\right\| _{q,p}^{\pi }=\left\{ 
\begin{array}{lll}
\left( \underset{E\in \pi }{\sum }\left\| f\chi _{_{E}}\right\|
_{q}^{p}\right) ^{\frac{1}{p}} & \text{si} & p<+\infty \\ 
\underset{E\in \pi }{\sup }\left\| f\chi _{_{E}}\right\| _{q} & \text{si} & 
p=+\infty%
\end{array}
\right. 
\end{equation*}

$\qquad \blacktriangleright \ L_{\left( q,p\right) }^{\pi }(G)=\left\{ f\in
L_{0}(G)\text{ / }\left\| f\right\| _{q,p}^{\pi }<+\infty \right\} $

Il est d\'{e}montr\'{e} dans $\cite{1}$ que :

$\qquad \qquad \bullet \left( L_{\left( q,p\right) }^{\pi }(G)\ ,\ \left\|
{}\right\| _{q,p}^{\pi }\right) $\ est un espace de Banach complexe.

$\qquad \qquad \bullet \ $si $1\leq q_{1}\leq q_{2}\leq p_{2}\leq p_{1}\leq
+\infty ,$ alors $L_{\left( q_{2},p_{2}\right) }^{\pi }(G)\subset L_{\left(
q_{1},p_{1}\right) }^{\pi }(G)$\bigskip

b) Consid\'{e}rons un voisinage ouvert, relativement compact et sym\'{e}%
trique $B$\ de $e$.

Posons :

$\qquad \blacktriangleright \ $Pour tout \'{e}l\'{e}ment $f$\ \ de $L_{0}(G)$%
,

\begin{equation*}
_{B}\left\| f\right\| _{q,p}=\left\{ 
\begin{array}{lll}
\left[ \int_{G}\left( \left\| f\chi _{_{yB}}\right\| _{q}\right)
^{p}d\lambda (y)\right] ^{\frac{1}{p}} & \text{si} & p<+\infty \\ 
\underset{y\in G}{\sup ess}\left\| f\chi _{_{yB}}\right\| _{q} & \text{si} & 
p=+\infty%
\end{array}
\right. \ ; 
\end{equation*}

$\qquad \blacktriangleright \ \left( L^{q},L^{p}\right) (G)=\left\{ f\in
L_{0}(G)\text{ / }_{B}\left\| f\right\| _{q,p}<+\infty \right\} .$

$\qquad \blacktriangleright \ $Pour tous \'{e}l\'{e}ments $f$ et $g$ de $%
L_{0}(G),$

\qquad \qquad\ $\left( f\ast g\right) (x)=\int_{G}f(y)g(y^{-1}x)d\lambda
(y) $\ \ en tout point $x$\ de $G$\ o\`{u} cela a un sens.

$\qquad \blacktriangleright \forall \left( a\ ,\ f\right) \in G\times
L_{0}(G),\ \left\{ 
\begin{array}{l}
f_{a}(x)=f(xa) \\ 
_{a}f(x)=f(a^{-1}x)%
\end{array}
\right. ,$

pour presque tout \'{e}l\'{e}ment $x$\ de $G$.\bigskip

Dans le paragraphe 2, nous d\'{e}montrons que pour tous \'{e}l\'{e}ments $p$
et $q$ de $\left[ 1\ ;\ +\infty \right] ,$

\qquad $\bullet \ $l'espace $\left( L^{q},L^{p}\right) \left( G\right) $ ne d%
\'{e}pend pas de $B\ $(voir proposition \ref{PROP: Equiv entre normes B( )
qp})$,$

\qquad $\bullet \ _{B}\left\| {}\right\| _{q,p}$ est une norme sur $%
L_{q,p}^{\pi }(G)$ \'{e}quivalente \`{a} la norme $\left\| {}\right\|
_{q,p}^{\pi }$ (voir proposition \ref{PROP: B( )qp<()piqp<cB()qp})

Ce deuxi\`{e}me r\'{e}sultat montre que l'espace $\left( L^{q},L^{p}\right)
\left( G\right) $ est identifiable \`{a} l'espace de Banach $L_{\left(
q,p\right) }^{\pi }(G)$.\bigskip

A la suite de F. Holland, de nombreux auteurs ont \'{e}tudi\'{e} ces espaces
dans le cas d'un groupe commutatif, en liaison avec les transformations de
Fourier.

c) Pour tout \'{e}l\'{e}ment $f$\ de $L_{0}(G),$ notons :\medskip

$\qquad \blacktriangleright \ \lambda _{f}$\ la fonction de distribution de $%
f$\ d\'{e}finie sur $\left[ 0\ ;\ +\infty \right[ $\ par :

\begin{equation*}
\lambda _{f}(s)=\lambda \left( \left\{ x\in G\text{ }/\text{ }\left| \text{ }%
f(x)\right| >s\right\} \right) , 
\end{equation*}

$\qquad \blacktriangleright \ f^{\ast }$\ sa fonction de r\'{e}arrangement d%
\'{e}croissante d\'{e}finie sur $\left[ 0\ ;\ +\infty \right[ $ par :\ 
\begin{equation*}
f^{\ast }(s)=\inf \left\{ \text{ }t>0\text{ }/\text{ }\lambda _{f}(t)\leq
s\right\} , 
\end{equation*}

$\qquad \blacktriangleright \ \left\| f\right\| _{q,p}^{\ast }=\left\{ 
\begin{array}{lll}
\left( \dfrac{p}{q}\int_{0}^{+\infty }\left( t^{\frac{1}{q}}f^{\ast
}(t)\right) ^{p}\dfrac{dt}{t}\right) ^{\frac{1}{p}} & \text{si} & 1\leq
q,p<+\infty \\ 
\underset{t>0}{\sup }t^{\frac{1}{q}}f^{\ast }(t) & \text{si} & 1\leq
q<p=+\infty \\ 
\underset{t>0}{\sup }f^{\ast }(t) & \text{si} & q=p=+\infty%
\end{array}
\right. ,\medskip $

$\qquad \blacktriangleright \ $Pour tous \'{e}l\'{e}ments $p$\ et $q$\ de $%
\left[ 1\ ;\ +\infty \right] ,\medskip $

\qquad $\qquad L^{q,p}(G)=\left\{ f\in L_{0}(G)\text{ / }\left\| f\right\|
_{q,p}^{\ast }<+\infty \right\} $ (\ espace de Lorentz ).\bigskip

d) Dans le cas o\`{u} $G$\ est un groupe homog\`{e}ne ( groupe simplement
connexe, nilpotent et muni d'une famille de dilatations $\left\{ \text{ }%
\delta _{r}\text{ },\text{ }r>0\text{ }\right\} $\ ) ,

notons :

$\qquad \blacktriangleright \ \left| \qquad \right| $\ une norme homog\`{e}%
ne fix\'{e}e dans $G$,

$\qquad \blacktriangleright \ \rho $\ sa dimension homog\`{e}ne,

$\qquad \blacktriangleright \gamma =\inf \left\{ c\in \text{ }\mathbb{R}%
_{+}^{\ast }\text{ }/\text{ }\left| xy\right| \leq c(\left| x\right| +\left|
y\right| \text{ }),\forall (x,y)\in G\times G\right\} .$

$\forall (x,r)\in G\times \mathbb{R}_{+}^{\ast },$\ 

$\qquad \blacktriangleright \ B_{\left( x,r\right) }=\left\{ \text{ }y\in G%
\text{ }/\text{ }\left| x^{-1}y\right| <r\text{ }\right\} $\ est la boule de
centre $x$ et de rayon $r$ ,

$\qquad \blacktriangleright \ \pi _{r}$\ d\'{e}signe une $B_{\left( e,\frac{r%
}{4\gamma ^{2}}\right) }-B_{\left( e,\frac{r}{4\gamma ^{2}}\right) }^{2}$\
partition uniforme de $G$.

$\forall \left( q,p,\alpha \right) \in \left[ 1\ ;\ +\infty \right] ^{3},$%
\begin{equation*}
\left\| f\right\| _{q,p,\alpha }=\underset{r>0}{\sup }\lambda \left(
B_{\left( e,r\right) }\right) ^{\frac{1}{\alpha }-\frac{1}{q}}\left\|
f\right\| _{q,p}^{\pi _{r}}\ ,\ \forall f\in L_{0}(G), 
\end{equation*}
et

\begin{equation*}
\left( L^{q},L^{p}\right) ^{\alpha }(G)=\left\{ \text{ }f\in L_{0}(G)\text{ }%
/\text{ }\left\| f\right\| _{q,p,\alpha }<+\infty \right\} . 
\end{equation*}

Nous remarquons que l'espace $\left( L^{q},L^{+\infty }\right) ^{\alpha }(G)$
co\"{\i}ncide avec l'espace classique de Morrey.

Dans le paragraphe 4 nous montrons entre autres choses que pour tout \'{e}l%
\'{e}ment $\left( q,p,\alpha \right) $ de $\left[ 1\ ;\ +\infty \right] ^{3}$
tel que $q\leq \alpha \leq p,$

$\qquad \blacktriangleright \ L^{\alpha }(G)\subset \left(
L^{q},L^{p}\right) ^{\alpha }(G)$ ( voir le corollaire \ref{cor: La(G) inclu
(Lq,Lp) a(G)}),

$\qquad \blacktriangleright \ \left( L^{q},L^{p}\right) ^{\alpha
}(G)=L^{\alpha }(G)$ si $\alpha =p$ ou $\alpha =q$ (voir les propositions %
\ref{PROP : L qqin=Lq}, \ref{PROP: Lqpq=Lq} et \ref{PROP: Lqpp=Lp}),

$\qquad \blacktriangleright \ L^{\alpha ,+\infty }(G)\subsetneq $ $\left(
L^{q},L^{p}\right) ^{\alpha }(G)$ si $q<\alpha <p$ (voir les propositions %
\ref{PROP: ( )qpa< ( )*a} et \ref{PROP : element de Lqpa}).

Ajoutons que dans le cas o\`{u} la condition $q\leq \alpha \leq p$ n'est pas
v\'{e}rifi\'{e}e, l'espace $\left( L^{q},L^{p}\right) ^{\alpha }(G)$ est 
\'{e}gal \`{a} $\left\{ O\right\} ,$ $O$ d\'{e}signant la classe des
fonctions nulles $\lambda -$presque partout (voir les propositions \ref%
{prop:cas a<q}, \ref{prop: p<a})\newpage

L'int\'{e}r\^{e}t des espaces $\left( L^{q},L^{p}\right) ^{\alpha }(G)$ est
entre autres mis en \'{e}vidence par le fait que les espaces $L^{\alpha }(G)$
de Lebesgue, $L^{\alpha ,+\infty }(G)$ de Lorentz, $\left( L^{q},L^{+\infty
}\right) ^{\alpha }(G)$ de Morrey en sont soit des sous-espaces, soit des
cas particuliers. Il faut aussi remarquer que de par leur d\'{e}finition,
ces espaces $\left( L^{q},L^{p}\right) ^{\alpha }(G),$ sont \'{e}troitement
li\'{e}s aux multiplicateurs et se pr\^{e}tent bien \`{a} l'\'{e}tude de l'op%
\'{e}rateur maximal fractionnaire de Hardy-Littlewood (voir \cite{2}). Nous
nous proposons dans un prochain travail d'\'{e}tudier les propri\'{e}t\'{e}s
de continuit\'{e} de cet op\'{e}rateur et de l'int\'{e}grale fractionnaire
dans ce cadre.

\section{ L'ESPACE DE BANACH $\left( L^{q},L^{p}\right) (G)$}

Rappelons quelques propositions et d\'{e}finitions tir\'{e}es de \cite{1}.

\begin{definition}
\label{DEF: Part-uni} Soient $U$\ et $V$\ deux voisinages ouverts\ et
relativement compacts de $e$ tels que $\overline{U}\subset V.$ On appelle $%
U-V$\ partition uniforme de $G,$ toute partition $\pi $\ de $G$\ en
boreliens telle que 
\begin{equation*}
\forall E\in \pi ,\exists x_{_{E}}\in E\ /\ x_{_{E}}U\subset E\subset
x_{_{E}}V. 
\end{equation*}
\bigskip
\end{definition}

\begin{proposition}
\label{PROP: Existence Part uniforme}Si $U$\ est un voisinage ouvert
relativement compact et sym\'{e}trique de $e,$\ alors il existe une $U-U^{2}$%
\ partition uniforme de $G$.
\end{proposition}

\bigskip

\begin{proposition}
\label{PROP: Equivalence normes uniformes}Soient $\pi $ et\ $\pi ^{\prime }$
deux\ partitions uniformes de $G$, $\ p$\ et $q$\ deux \'{e}l\'{e}ments de $%
\left[ 1\ ;\ +\infty \right] $, il existe un r\'{e}el positif $M=M(\pi ,\pi
^{\prime })$\ tel que pour tout \'{e}l\'{e}ment $f$\ de $L_{0}(G)$, nous
ayons 
\begin{equation*}
\left\| f\right\| _{q,p}^{\pi ^{\prime }}\leq M\left\| f\right\| _{q,p}^{\pi
}. 
\end{equation*}
\end{proposition}

\bigskip

\begin{proposition}
\label{PROP: Toute translaté rencontreN(K,L)}Soient $\pi $\ une $U-V$\
partition uniforme de $G$, $K$\ et $L$\ deux bor\'{e}liens de $G$\
relativement compacts. Il existe un r\'{e}el $n_{\pi }(K,L)\geq 1$\ tel que
toute translat\'{e}e \`{a} gauche de $L$\ rencontre au plus $n_{\pi }(K,L)$\
ensembles $x_{E}K$, avec $E$\ \'{e}l\'{e}ment de $\pi $. Nous pouvons
prendre 
\begin{equation}
n_{\pi }(K,L)=\dfrac{\lambda \left( LK^{-1}U\right) }{\lambda \left(
U\right) }.  \label{2.1}
\end{equation}
\end{proposition}

\bigskip

\begin{proposition}
\label{PROP: n(V,K) element pie}Soient $\pi $\ une $U-V$\ partition uniforme
de $G$ et $K$\ un bor\'{e}lien de $G$\ relativement compact. Toute translat%
\'{e}e \`{a} gauche de $K$\ rencontre au plus $n_{\pi }(V,K)$\ \'{e}l\'{e}%
ments de $\pi $.
\end{proposition}

\bigskip

Enon\c{c}ons aussi ce lemme dont la preuve est imm\'{e}diate.\bigskip

\begin{lemma}
\label{Lem:Incl}Soient $N,K$ et $O$ des sous-ensembles de $G$.

Si $\lambda \left( N\right) =0,\ K$ est relativement compact et $O$ est
ouvert, alors il existe un sous-ensemble fini $\left\{ y_{1};y_{2};\ldots
;y_{n}\right\} $ de $K\backslash N$ tel que $K\subset \underset{i=1}{\overset%
{n}{\cup }}y_{i}O.$ 
\end{lemma}

\bigskip

A l'aide de ces r\'{e}sultats, nous allons d\'{e}montrer l'\'{e}quivalence
des normes $\left\| {}\right\| _{q,p}^{\pi }$\ et $_{B}\left\| {}\right\|
_{q,p},$ avec $B$ un voisinage ouvert relativement compact et sym\'{e}trique
de $e$.\ Mais pour cela, nous avons besoin des propositions suivantes :
\bigskip

\begin{proposition}
\label{PROP: continuite de fconvqkiE}Soient $(p,q)$\ un \'{e}l\'{e}ment de $%
\left[ 1\ ;\ +\infty \right] ^{2}$ et $B$ un voisinage ouvert relativement
compact et sym\'{e}trique de $e$. Pour tout \'{e}l\'{e}ment $f$ de $\left(
L^{q},L^{p}\right) (G)$,

a) $\left( \left| f\right| ^{q}\ast \chi _{_{B}}\right) (x)$\ est fini pour
tout $x$\ \'{e}l\'{e}ment de $G.$

b) $\left| f\right| ^{q}\ast \chi _{_{B}}$\ est continue sur $G$.
\end{proposition}

\begin{proof}
\ Soient $f$ un \'{e}l\'{e}ment de $\left( L^{q},L^{p}\right) (G)$\ et $x$\
un \'{e}l\'{e}ment de $G.$

\textbf{a)} Il existe un sous-ensemble $\lambda -$mesurable $N$ de $G$ tel
que : 
\begin{equation*}
\left\{ 
\begin{array}{l}
\lambda \left( N\right) =0 \\ 
\forall y\in G\backslash N,\ \left\| f\chi _{_{yB}}\right\| _{q}<+\infty%
\end{array}
\right. . 
\end{equation*}

D'apr\`{e}s le lemme \ref{Lem:Incl}, il existe un sous-ensemble fini $%
\left\{ y_{1};y_{2};\ldots ;y_{n}\right\} $ de $xB\backslash N$ tel que $%
xB\subset \cup y_{i}B.$

Par cons\'{e}quent, 
\begin{equation*}
\left( \left| f\right| ^{q}\ast \chi _{_{B}}\right) (x)\leq \underset{i=1}{%
\overset{n}{\sum }}\left\| f\chi _{_{y_{i}B}}\right\| _{q}^{q}<+\infty . 
\end{equation*}

\textbf{b)} On consid\`{e}re une suite $\left( x_{n}\right) _{n\geq 0\text{ }%
}$d'\'{e}l\'{e}ments de $G$\ qui converge vers $x.$

\ Il existe un voisinage compact $K$\ de $e$, tel que 
\begin{equation*}
x_{n}x^{-1}\in K{\large \ }\forall n\in \mathbb{N}{\large .} 
\end{equation*}

\ Par cons\'{e}quent,

$\underset{n\rightarrow +\infty }{\lim }\left| \left( \left| f\right|
^{q}\ast \chi _{_{B}}\right) (x_{n})-\left( \left| f\right| ^{q}\ast \chi
_{_{B}}\right) (x)\right| \leq \underset{n\rightarrow +\infty }{\lim }%
\left\| _{xx_{n}^{-1}}\left( \left| f\right| ^{q}\chi _{_{Kx\overline{B}%
}}\right) -\left( \left| f\right| ^{q}\chi _{_{Kx\overline{B}}}\right)
\right\| _{1}=0.$
\end{proof}

\begin{proposition}
\label{PROP: Equiv entre normes B( ) qp}Soient $B_{1}$\ et $B_{2}$\ deux
voisinages ouverts, sym\'{e}triques et relativement compacts de $e$, $p$\ et 
$q$\ deux \'{e}l\'{e}ments de $\left[ 1\ ;\ +\infty \right] $. Il existe une
constante r\'{e}elle $C=C(B_{1},B_{2})$\ telle que :\ 
\begin{equation*}
_{B_{2}}\left\| f\right\| _{q,p}\leq C\text{ }_{B_{1}}\left\| f\right\|
_{q,p}, 
\end{equation*}
pour tout \'{e}l\'{e}ment $f$ de $L_{0}(G).$
\end{proposition}

\begin{proof}
Soit $f$\ un \'{e}l\'{e}ment de $L_{0}(G)$. Il existe un sous-ensemble fini $%
\left\{ y_{1};y_{2};\ldots ;y_{n}\right\} $ de $B_{2}$ tel que $B_{2}\subset 
\underset{i=1}{\overset{n}{\cup }}B_{1}y_{i}.$

Par suite, pour tout \'{e}l\'{e}ment $y$ de $G,$ 
\begin{equation*}
\left\{ 
\begin{array}{lll}
\left( \left| f\right| ^{q}\ast \chi _{_{B_{2}}}\right) ^{\frac{1}{q}%
}(y)\leq \underset{i=1}{\overset{n}{\sum }}\left[ \left( \left| f\right|
^{q}\ast \chi _{_{B_{1}}}\right) (yy_{i}^{-1})\right] ^{\frac{1}{q}} & \text{%
si} & q<+\infty \\ 
\left\| f\chi _{_{yB_{2}}}\right\| _{+\infty }\leq \underset{i=1}{\overset{n}%
{\sum }}\left\| f\chi _{_{yy_{i}B_{1}}}\right\| _{+\infty } & \text{si} & 
q=+\infty%
\end{array}
.\right. 
\end{equation*}
En prenant la norme $L^{p}(G)$\ des deux \ membres de chacune des in\'{e}%
galit\'{e}s, nous obtenons le r\'{e}sultat d\'{e}sir\'{e}.
\end{proof}

\bigskip

On d\'{e}duit de ce qui pr\'{e}c\`{e}de que les normes $_{B}\left\|
{}\right\| _{q,p}$ sont \'{e}quivalentes pour diff\'{e}rents choix de $B$.%
\bigskip

\begin{remark}
Soit $p$\ un \'{e}l\'{e}ment de l'intervalle $\left[ 1\ ;\ +\infty \right] $%
, on a $\left( L^{p},L^{p}\right) (G)=L^{p}(G)$.
\end{remark}

\begin{proof}
Elle est imm\'{e}diate \`{a} partir de la d\'{e}finition de $_{B}\left\|
\quad \right\| _{q,p}.$
\end{proof}

\bigskip

\begin{proposition}
\label{PROP: B( )qp<()piqp<cB()qp}Soient $p$\ et $q$\ deux \'{e}l\'{e}ments
de $\left[ 1\ ;\ +\infty \right] .$ Pour tout voisinage ouvert relativement
compact et sym\'{e}trique $B$ de $e$ et pour toute partition uniforme $\pi $
de $G$, il existe deux constantes r\'{e}elles $C_{1}=C_{1}(B,\pi )$\ et $%
C_{2}=C_{2}(B,\pi )$\ telles que pour tout \'{e}l\'{e}ment $f$ de $L_{0}(G)$%
,\ 
\begin{equation*}
C_{1}\text{ }_{B}\left\| f\right\| _{q,p}\leq \left\| f\right\| _{q,p}^{\pi
}\leq C_{2}\text{ }_{B}\left\| f\right\| _{q,p}. 
\end{equation*}
\end{proposition}

\begin{proof}
\ Soient $B_{1}$\ et $B_{2}$\ deux voisinages ouverts sym\'{e}triques et
relativement compacts de $e$, tels que 
\begin{equation*}
B_{1}^{2}\subset B_{2}\ \text{et}\ B_{2}^{2}\subset B. 
\end{equation*}

D\'{e}signons par $\pi $\ une $B_{1}-B_{1}^{2}$\ partition uniforme de $G$.
Pour tout \'{e}l\'{e}ment $E$ de $\pi ,$\ 
\begin{equation*}
\left\{ 
\begin{array}{l}
\exists x_{E}^{0}\in E\text{ }/\text{ }x_{E}^{0}B_{1}\subset E\subset
x_{E}^{0}B_{1}^{2} \\ 
\forall x\in E,\text{ }E\subset xB%
\end{array}
\right. . 
\end{equation*}

Posons : 
\begin{equation*}
T(E)=\left\{ E^{\prime }\in \pi \text{ / }E\cap E^{\prime }B\neq \emptyset
\right\} \ \text{et}\ T_{y}=\left\{ E\in \pi \text{ / }E\cap yB\neq
\emptyset \right\} , 
\end{equation*}
pour tous \'{e}l\'{e}ments $E$\ de $\pi $, et $y$\ de $G.$

D'apr\`{e}s la proposition \ref{PROP: n(V,K) element pie},\ 

$card\ T(E)\leq n_{\pi }(B_{1}^{2}\ ;\ B_{1}^{2}B),$ et $\ card\ T_{y}\leq
n_{\pi }(B_{1}^{2},B).$

Soit $f$\ un \'{e}l\'{e}ment de $L_{0}(G)$. Nous avons :

D'une part 
\begin{equation*}
\left\{ 
\begin{array}{lllll}
\left\| f\chi _{_{E}}\right\| _{q}\leq \left[ \lambda \left( B_{1}\right)
^{-1}\int_{E}\left( \left| f\right| ^{q}\ast \chi _{_{B}}\right) ^{\frac{p}{%
q}}(x)d\lambda (x)\right] ^{\frac{1}{p}} & \text{si} & p<+\infty & \text{et}
& q<+\infty \\ 
\left\| f\chi _{_{E}}\right\| _{q}\leq \left\| \left( \left| f\right|
^{q}\ast \chi _{_{B}}\right) ^{\frac{1}{q}}\right\| _{+\infty } & \text{si}
& p=+\infty & \text{et} & q<+\infty \\ 
\left\| f\chi _{_{E}}\right\| _{+\infty }\leq \left[ \lambda \left(
B_{1}\right) ^{-1}\int_{E}\left\| f\chi _{_{yB}}\right\| _{+\infty
}^{p}d\lambda (y)\right] ^{\frac{1}{p}} & \text{si} & \text{\ }p<+\infty & 
\text{et} & q=+\infty%
\end{array}
\right. , 
\end{equation*}
pour tout \'{e}l\'{e}ment $E$ de $\pi .$\ 

D'o\`{u} en prenant la norme $\ell ^{p}\left( \pi \right) $\ des deux
membres de chacune des in\'{e}galit\'{e}s, nous obtenons

\begin{equation*}
\left\| f\right\| _{q,p}^{\pi }\leq C_{2}\ _{B}\left\| f\right\| _{q,p},%
\text{ o\`{u} }C_{2}\ \text{est une constante ind\'{e}pendante de }f\text{.} 
\end{equation*}

D'autre part,\ 
\begin{equation*}
\left\{ 
\begin{array}{lllll}
\left( \left| f\right| ^{q}\ast \chi _{_{B}}\right) ^{\frac{1}{q}}(y)\leq
n_{\pi }(B_{1}^{2},B)^{\frac{1}{q}}\left( \underset{E_{1}\in \pi }{\sum }%
\left\| f\chi _{_{E_{1}\cap yB}}\right\| _{q}^{p}\right) ^{\frac{1}{p}} & 
\text{si} & q<+\infty & \text{et} & p<+\infty \\ 
\left( \left| f\right| ^{q}\ast \chi _{_{B}}\right) ^{\frac{1}{q}}(y)\leq
n_{\pi }(B_{1}^{2},B_{1}^{2}B)\left\| f\right\| _{q,+\infty }^{\pi } & \text{%
si} & q<+\infty & \text{et} & p=+\infty \\ 
\left\| f\chi _{_{yB}}\right\| _{+\infty }\leq \underset{E\in T_{y}}{\sum }%
\left\| f\chi _{_{E}}\right\| _{+\infty }\chi _{_{x_{E}^{0}B_{1}^{2}B}}(y) & 
\text{si} & q=+\infty & \text{et} & p=+\infty%
\end{array}
\right. , 
\end{equation*}
pour $\lambda -$presque tout $y$ dans $G$.

D'o\`{u} en prenant la norme $L^{p}(G)$\ des deux membres de chacune des in%
\'{e}galit\'{e}s,

nous obtenons$\qquad _{B}\left\| f\right\| _{q,p}\leq C_{1}^{-1}\left\|
f\right\| _{q,p}^{\pi }$
\end{proof}

\bigskip

\begin{remark}
Soient $p$ et $q$ deux \'{e}l\'{e}ments de $\left[ 1\ ;\ +\infty \right] .$

1- La proposition \ref{PROP: B( )qp<()piqp<cB()qp} montre que l'espace $%
\left( L^{q},L^{p}\right) (G)$ est ind\'{e}pendante du voisinage ouvert sym%
\'{e}trique et relativement compact $B$ utilis\'{e} pour sa d\'{e}finition
(ce qui justifie notre notation)

2- Soit $\pi $ une partition uniforme de $G$. Sachant que $\left(
L_{q,p}^{\pi }(G),\left\| \quad \right\| _{q,p}^{\pi }\right) $ est un
espace de Banach complexe, il est ais\'{e} de voir \`{a} partir de la
proposition \ref{PROP: B( )qp<()piqp<cB()qp} que :

\qquad pour tout voisinage ouvert sym\'{e}trique et relativement compact $B$
de $e$, $_{B}\left\| \quad \right\| _{q,p}$ est une norme sur $\left(
L^{q},L^{p}\right) (G)=L_{q,p}^{\pi }(G)$, \'{e}quivalente \`{a} la norme $%
\left\| \quad \right\| _{q,p}^{\pi }.$
\end{remark}

\section{ ESPACE $\left( L^{q},L^{p}\right) ^{\protect\alpha }(G)$}

Dans tout ce paragraphe, $G$\ est un groupe homog\`{e}ne.

Au vu de la d\'{e}finition de l'espace $\left( L^{q},L^{p}\right) ^{\alpha
}(G),$ la proposition \ref{PROP: B( )qp<()piqp<cB()qp} prend la forme plus pr%
\'{e}cise suivante :

\bigskip

\begin{proposition}
\label{PROP: équivalence BusbyetB( )qp}Soient $(p,q)$\ un \'{e}l\'{e}ment de 
$\left[ 1\ ;\ +\infty \right] ^{2}$, $r$\ un r\'{e}el strictement positif et 
$f$\ un \'{e}l\'{e}ment de $\left( L^{q},L^{p}\right) (G).$ Alors il existe
deux constantes $C_{1}=C_{1}\left( \gamma ,\rho \right) $ et $%
C_{2}=C_{2}\left( \gamma ,\rho \right) $ telles que\ nous ayons : 
\begin{equation*}
C_{1}\left\| f\right\| _{q,p}^{\pi _{r}}\leq r^{-\frac{\rho }{p}}\text{ }%
_{B_{\left( e,r\right) }}\left\| f\right\| _{q,p}\leq C_{2}\left\| f\right\|
_{q,p}^{\pi _{r}}. 
\end{equation*}
\end{proposition}

\begin{proof}
Posons 
\begin{equation*}
B_{1}=B_{\left( e,\frac{r}{4\gamma ^{2}}\right) },B_{2}=B_{\left( e,\frac{r}{%
4\gamma ^{2}}\right) }^{2}\text{ et }B=B_{\left( e,r\right) }. 
\end{equation*}
Remarquons que : 
\begin{equation*}
B_{1}^{2}B\left( B_{1}^{2}\right) ^{-1}B_{1}\subset B_{\left( e,r(\gamma
^{3}+\frac{3}{4}\gamma +\frac{1}{2})\right) }\text{ et }B\left(
B_{1}^{2}\right) ^{-1}B_{1}\subset B_{\left( e,r(\gamma ^{2}+\frac{3}{4}%
)\right) }. 
\end{equation*}

D'apr\`{e}s la relation (\ref{2.1}) de la proposition \ref{PROP: Toute
translaté rencontreN(K,L)} et les propri\'{e}t\'{e}s de la mesure $\lambda $%
,\ nous avons : 
\begin{equation*}
n_{\pi _{r}}\left( B_{1}^{2},B_{1}^{2}B\right) <\dfrac{\lambda \left(
B_{\left( e,r(\gamma ^{3}+\frac{3}{4}\gamma +\frac{1}{2})\right) }\right) }{%
\lambda \left( B_{\left( e,\frac{r}{4\gamma ^{2}}\right) }\right) }=\dfrac{%
\left( r(\gamma ^{3}+\frac{3}{4}\gamma +\frac{1}{2})\right) ^{\rho }}{\left( 
\frac{r}{4\gamma ^{2}}\right) ^{\rho }}=\left( 4\gamma ^{5}+3\gamma
^{3}+2\gamma ^{2}\right) ^{\rho }, 
\end{equation*}
\begin{equation*}
n_{\pi _{r}}\left( B_{1}^{2},B\right) <\dfrac{\lambda \left( B_{\left(
e,r(\gamma ^{2}+\frac{3}{4})\right) }\right) }{\lambda \left( B_{\left( e,%
\frac{r}{4\gamma ^{2}}\right) }\right) }=\left( 4\gamma ^{4}+3\gamma
^{2}\right) ^{\rho }, 
\end{equation*}
\begin{equation*}
\lambda \left( B\right) =r^{\rho },\ \lambda \left( B_{1}\right) =\left( 
\frac{r}{4\gamma ^{2}}\right) ^{\rho }\ \text{et}\ \lambda \left(
B_{2}\right) =\left( \frac{r}{2\gamma }\right) ^{\rho }. 
\end{equation*}

Le r\'{e}sultat annonc\'{e} s'obtient en reprenant la d\'{e}monstration de
la proposition \ref{PROP: B( )qp<()piqp<cB()qp} et en tenant compte des pr%
\'{e}cisions ci-dessus.
\end{proof}

\bigskip

\begin{proposition}
Soit $\left( \alpha ,p,q\right) $\ un \'{e}l\'{e}ment de $\left[ 1\ ;\
+\infty \right] ^{3}$.

\ a) $\left( L^{q},L^{p}\right) ^{\alpha }(G)$, est un sous-espace vectoriel
complexe de $\left( L^{q},L^{p}\right) (G)$.

\ b) L'application $f\longmapsto \left\| f\right\| _{q,p,\alpha }$\ d\'{e}%
finit une norme sur $\left( L^{q},L^{p}\right) ^{\alpha }(G)$.
\end{proposition}

\begin{proof}
\ Elle est imm\'{e}diate \`{a} partir des d\'{e}finitions des normes $%
\left\| {}\right\| _{q,p}^{\pi _{r}}$\ et $\left\| {}\right\| _{q,p,\alpha
}. $
\end{proof}

\begin{proposition}
Pour tout triplet $\left( \alpha ,p,q\right) $\ d'\'{e}l\'{e}ments de $\left[
1\ ;\ +\infty \right] .$

$\left( L^{q},L^{p}\right) ^{\alpha }(G)$\ muni de la norme $\left\|
{}\right\| _{q,p,\alpha },\ $est un espace de Banach complexe.
\end{proposition}

\begin{proof}
Soit $\left( f_{n}\right) _{n\in \mathbb{N}},$\ une suite de Cauchy dans $%
\left( L^{q},L^{p}\right) ^{\alpha }(G).$ Pour tous entiers naturels $n$\ et 
$m$, nous avons : 
\begin{equation*}
\left\| f_{n}-f_{m}\right\| _{q,p}^{\pi _{1}}\leq \left\|
f_{n}-f_{m}\right\| _{q,p,\alpha }. 
\end{equation*}

Par cons\'{e}quent, $\left( f_{n}\right) _{n\in \mathbb{N}}$\ est une suite
de Cauchy dans l'espace de Banach $\left( L^{q},L^{p}\right) (G)$, et par
suite y converge vers un \'{e}l\'{e}ment $f$.

De m\^{e}me, $\left( \left\| f_{n}\right\| _{q,p,\alpha }\right) _{n\in 
\mathbb{N}}$ \'{e}tant de Cauchy dans $\mathbb{R},$ y converge vers un r\'{e}%
el $M\geq 0 $.

Pour tout r\'{e}el $r>0$\ et tout entier $m\geq 0$, 
\begin{equation*}
\begin{array}{lll}
\lambda \left( B_{\left( e,r\right) }\right) ^{\frac{1}{\alpha }-\frac{1}{q}%
}\left\| f\right\| _{q,p}^{\pi _{r}} & = & \lambda \left( B_{\left(
e,r\right) }\right) ^{\frac{1}{\alpha }-\frac{1}{q}}\left\|
f_{m}-f_{m}+f\right\| _{q,p}^{\pi _{r}} \\ 
& \leq & \lambda \left( B_{\left( e,r\right) }\right) ^{\frac{1}{\alpha }-%
\frac{1}{q}}\left\| f_{m}\right\| _{q,p}^{\pi _{r}}+\lambda \left( B_{\left(
e,r\right) }\right) ^{\frac{1}{\alpha }-\frac{1}{q}}\left\| f_{m}-f\right\|
_{q,p}^{\pi _{r}} \\ 
& \leq & \left\| f_{m}\right\| _{q,p,\alpha }+\lambda \left( B_{\left(
e,r\right) }\right) ^{\frac{1}{\alpha }-\frac{1}{q}}\left\| f_{m}-f\right\|
_{q,p}^{\pi _{r}}.%
\end{array}
\end{equation*}

\ Donc pour tout r\'{e}el $r>0$, \ $\lambda \left( B_{\left( e,r\right)
}\right) ^{\frac{1}{\alpha }-\frac{1}{q}}\left\| f\right\| _{q,p}^{\pi
_{r}}\leq M.$

Ainsi, $f$\ est un \'{e}l\'{e}ment de $\left( L^{q},L^{p}\right) ^{\alpha
}(G),$ puisque $\underset{r>0}{\sup }\lambda \left( B_{\left( e,r\right)
}\right) ^{\frac{1}{\alpha }-\frac{1}{q}}\left\| f\right\| _{q,p}^{\pi
_{r}}<+\infty .$

Soit un r\'{e}el $\varepsilon >0$.

Il existe un entier naturel $m_{0}$\ tel que pour tous entiers naturels $%
m^{\prime }>m_{0}$\ et $m">m_{0},$\ nous ayons : 
\begin{equation*}
\lambda \left( B_{\left( e,r\right) }\right) ^{\frac{1}{\alpha }-\frac{1}{q}%
}\left\| f_{m^{\prime }}-f_{m"}\right\| _{q,p}^{\pi _{r}}\leq \left\|
f_{m\prime }-f_{m"}\right\| _{q,p,\alpha }<\varepsilon ,\ \forall r\in 
\mathbb{R}_{+}^{\ast }. 
\end{equation*}

Donc, pour tout r\'{e}el $r>0$\ et tout $m^{\prime }>m_{0}$,\ 
\begin{equation*}
\lambda \left( B_{\left( e,r\right) }\right) ^{\frac{1}{\alpha }-\frac{1}{q}%
}\left\| f_{m^{\prime }}-f\right\| _{q,p}^{\pi _{r}}=\underset{m"\rightarrow
+\infty }{\lim }\lambda \left( B_{\left( e,r\right) }\right) ^{\frac{1}{%
\alpha }-\frac{1}{q}}\left\| f_{m^{\prime }}-f_{m"}\right\| _{q,p}^{\pi
_{r}}<\varepsilon . 
\end{equation*}

Par cons\'{e}quent, pour tout $m^{\prime }>m_{0}$,\ $\left\| f_{m^{\prime
}}-f\right\| _{q,p,\alpha }<\varepsilon .$\ 

D'o\`{u} $\left( f_{n}\right) _{n\in \mathbb{N}}$\ converge dans $\left(
L^{q},L^{p}\right) ^{\alpha }(G)$\ vers $f.$
\end{proof}

\bigskip

\begin{proposition}
Soit $\left( p,q,\alpha \right) $\ un \'{e}l\'{e}ment de $\left[ 1\ ;\
+\infty \right] ^{3}$ avec $1\leq q\leq \alpha \leq p\leq +\infty $. Il
existe deux constantes $C_{1}$ et $C_{2}$\ telles que pour tout \'{e}l\'{e}%
ment $f$ de $\left( L^{q},L^{p}\right) (G)$, nous ayons :\ 
\begin{equation*}
C_{1}\left\| f\right\| _{q,p,\alpha }\leq \underset{r>0}{\sup }\lambda
\left( B_{\left( e,r\right) }\right) ^{\frac{1}{\alpha }-\frac{1}{q}-\frac{1%
}{p}}\text{ }_{B_{\left( e,r\right) }}\left\| f\right\| _{q,p}\leq
C_{2}\left\| f\right\| _{q,p,\alpha }. 
\end{equation*}
\end{proposition}

\begin{proof}
\ \ Elle est imm\'{e}diate \`{a} partir de la proposition \ \ref{PROP:
équivalence BusbyetB( )qp}.
\end{proof}

\bigskip

\begin{proposition}
\label{prop:cas a<q}Soient $\alpha ,q$\ et $p$\ des \'{e}l\'{e}ments de $%
\left[ 1\ ;\ +\infty \right] $.

Si $\alpha <q$, alors $\left( L^{q},L^{p}\right) ^{\alpha }(G)=\left\{
O\right\} $, o\`{u} $O$\ d\'{e}signe la classe des fonctions nulles $\lambda
-$presque partout sur $G$.
\end{proposition}

\begin{proof}
Soit $f$\ un \'{e}l\'{e}ment de $\left( L^{q},L^{p}\right) ^{\alpha }(G)$.

Posons$\qquad \ A=\left\| f\right\| _{q,p,\alpha }$.

Pour tout r\'{e}el strictement positif $r$,\ 
\begin{equation*}
\left\| f\right\| _{q,p}^{\pi _{r}}\leq \frac{A}{\lambda \left( B_{\left(
e,r\right) }\right) ^{\frac{1}{\alpha }-\frac{1}{q}}}\ ; 
\end{equation*}
et puisque$\ \frac{1}{\alpha }-\frac{1}{q}>0$,\ 
\begin{equation*}
\underset{r\rightarrow +\infty }{\lim }\left\| f\right\| _{q,p}^{\pi
_{r}}=0. 
\end{equation*}

Si\textbf{\ } $p<q$, alors nous avons pour tout r\'{e}el strictement positif 
$r,$%
\begin{equation*}
\left\| f\right\| _{q}=\left\| f\right\| _{q,q}^{\pi _{r}}\leq \left\|
f\right\| _{q,p}^{\pi _{r}}. 
\end{equation*}
Donc $f=O.$

Supposons $q<p$.

Puisque $f$\ est un \'{e}l\'{e}ment de $\left( L^{q},L^{p}\right) ^{\alpha
}(G)$, $\ f$\ est localement int\'{e}grable. Consid\'{e}rons un
sous-ensemble compact $K$\ de $G$.

Posons\ : 
\begin{equation*}
r_{0}=\inf \left\{ r>0\text{ }/\text{ }K\subset B_{\left( e,r\right)
}\right\} , 
\end{equation*}
et pour tout r\'{e}el r, 
\begin{equation*}
T_{K,r}=\left\{ E\in \pi _{r}\text{ }/\text{ }E\cap K\neq 0\right\} . 
\end{equation*}

Pour tout r\'{e}el $r>r_{0}$, 
\begin{equation*}
K\subset B_{\left( e,r\right) }, 
\end{equation*}
et 
\begin{equation*}
\left\| f\chi _{_{K}}\right\| _{q}^{q}=\underset{E\in T_{K,r}}{\sum }\left\|
f\chi _{_{E}}\right\| _{q}^{q}\leq \left( 4\gamma ^{4}+3\gamma ^{2}\right) ^{%
\frac{\rho q}{p}}\left( \underset{E\in T_{K,r}}{\sum }\left\| f\chi
_{_{E}}\right\| _{q}^{p}\right) ^{\frac{q}{p}}\leq \left( 4\gamma
^{4}+3\gamma ^{2}\right) ^{\frac{\rho q}{p}}\left( \left\| f\right\|
_{q,p}^{\pi _{r}}\right) ^{q}. 
\end{equation*}

\ Donc,\ 
\begin{equation*}
\left\| f\chi _{_{K}}\right\| _{q}\leq \left( 4\gamma ^{4}+3\gamma
^{2}\right) ^{\frac{\rho }{p}}\left( \underset{r\rightarrow +\infty }{\lim }%
\left\| f\right\| _{q,p}^{\pi _{r}}\right) =0 
\end{equation*}

Et puisque $G$\ est $\sigma $-compact, $f$\ est nulle $\lambda -$presque
partout sur $G$.

Le cas o\`{u} $p=q$ est trivial, car $\left\| f\right\| _{p,p}^{\pi
}=\left\| f\right\| _{p}.\bigskip $
\end{proof}

\begin{lemma}
Supposons que $1\leq p<\theta <+\infty $ et posons : 
\begin{equation*}
\forall y\in G\ h(y)=\left| y\right| ^{-\frac{\rho }{\theta ^{\prime }}} 
\end{equation*}
Alors, pour tout \'{e}l\'{e}ment a de $\mathbb{R}_{+}^{\ast },$%
\begin{equation*}
k(\theta ;p;a)=\left\| h\chi _{_{G\backslash B_{\left( e,a\right)
}}}\right\| _{p^{\prime }}<+\infty . 
\end{equation*}
\end{lemma}

\begin{proof}
a) Supposons que $1<p.$%
\begin{equation*}
1<p<\theta <+\infty \Longrightarrow 1<\theta ^{\prime }<p^{\prime }<+\infty
\Longrightarrow \rho <p^{\prime }\frac{\rho }{\theta ^{\prime }}<+\infty . 
\end{equation*}
Par cons\'{e}quent, d'apr\`{e}s le corollaire 1.6 de \cite{4}, il existe un
nombre r\'{e}el $C_{0}$ d\'{e}pendant uniquement de $G$, tel que : 
\begin{equation*}
\forall b\in \left] a\ ;\ +\infty \right[ \quad \int_{B_{\left( e,b\right)
}\backslash B_{\left( e,a\right) }}\left| y\right| ^{-p^{\prime }\frac{\rho 
}{\theta ^{\prime }}}d\lambda (y)=\dfrac{C_{0}}{d}\left(
a^{-d}-b^{-d}\right) , 
\end{equation*}
avec $d=\rho \left( \dfrac{p^{\prime }}{\theta ^{\prime }}-1\right) .$

En faisant tendre $b$ vers +$\infty ,$ il vient 
\begin{equation*}
\left\| h\chi _{_{G\backslash B_{\left( e,a\right) }}}\right\| _{p^{\prime
}}=\left( \frac{C_{0}a^{-d}}{d}\right) ^{\frac{1}{p^{\prime }}}<+\infty . 
\end{equation*}

b) Si $p=1$ alors $p^{\prime }=+\infty $ et

$\qquad \qquad \qquad \qquad \qquad \qquad \qquad \left\| h\chi
_{_{G\backslash B_{\left( e,a\right) }}}\right\| _{+\infty }=a^{-\frac{\rho 
}{\theta ^{\prime }}}<+\infty .\bigskip $
\end{proof}

\begin{lemma}
\label{Lem: B()q'p'<kr}Supposons $1\leq p<\theta <+\infty ,\ 1\leq q\leq
+\infty $ et posons 
\begin{equation*}
\forall x\in G\quad g(x)=\left( \gamma +\gamma \left| x\right| \right) ^{-%
\frac{\rho }{\theta ^{\prime }}}. 
\end{equation*}
Alors pour tout \'{e}l\'{e}ment $\varepsilon $ de $\left] 0\ ;\ 2\gamma %
\right[ ,$ nous avons : 
\begin{equation*}
\forall r\in \left] 0\ ;\ \frac{\varepsilon }{2\gamma }\right[ \quad
_{B_{\left( e,r\right) }}\left\| g\chi _{_{G\backslash B_{\left(
e,\varepsilon \right) }}}\right\| _{q^{\prime },p^{\prime }}\leq k(\theta \
;\ p\ ;\ \frac{\varepsilon }{2\gamma })\lambda \left( B_{\left( e,r\right)
}\right) ^{\frac{1}{q^{\prime }}}. 
\end{equation*}
\end{lemma}

\begin{proof}
Consid\'{e}rons un \'{e}l\'{e}ment $\left( \varepsilon ,r\right) $ de $%
\mathbb{R}_{+}^{\ast }\times \mathbb{R}_{+}^{\ast }$ tel que $\varepsilon
<2\gamma $ et $r<\frac{\varepsilon }{2\gamma }.$

1$^{er}cas$ $q=1$ ; c'est-\`{a}-dire que $q^{\prime }=+\infty .$

Consid\'{e}rons deux \'{e}l\'{e}ments $x$ et $y$ de $G$ v\'{e}rifiant 
\begin{equation*}
\chi _{_{G\backslash B_{\left( e,\varepsilon \right) }}}\left( x\right) \chi
_{_{yB_{\left( e,r\right) }}}\left( x\right) \neq 0 
\end{equation*}
Nous avons

$y^{-1}x\in B_{\left( e,r\right) }$ et donc $\left| y^{-1}x\right| <r$

$x\notin B_{\left( e,\varepsilon \right) }$ et donc $\varepsilon <\left|
x\right| \leq \gamma \left( \left| y\right| +\left| x\right| \right) <\gamma
\left| y\right| +\gamma r.$

$\varepsilon <\gamma \left| y\right| +\gamma r\Longrightarrow \dfrac{%
\varepsilon }{\gamma }-r<\left| y\right| \Longrightarrow \dfrac{\varepsilon 
}{2\gamma }<\left| y\right| $

$\left| y\right| \leq \gamma \left( \left| x\right| +\left| x^{-1}y\right|
\right) <\gamma \left| x\right| +\gamma r<\gamma \left| x\right| +\frac{%
\varepsilon }{2}<\gamma \left| x\right| +\gamma \Longrightarrow \left(
\gamma \left| x\right| +\gamma \right) ^{-\frac{\rho }{\theta ^{\prime }}%
}<\left| y\right| ^{-\frac{\rho }{\theta ^{\prime }}}.$

Ainsi 
\begin{equation*}
\left\| g\chi _{_{G\backslash B_{\left( e,\varepsilon \right) }}}\chi
_{_{yB_{\left( e,r\right) }}}\right\| _{+\infty }\leq \left| y\right| ^{-%
\frac{\rho }{\theta ^{\prime }}}\chi _{_{G\backslash B_{\left( e,\frac{%
\varepsilon }{2\gamma }\right) }}}(y), 
\end{equation*}
et par suite

$_{B_{\left( e,r\right) }}\left\| g\chi _{_{G\backslash B_{\left(
e,\varepsilon \right) }}}\right\| _{+\infty ,p^{\prime }}\leq k(\theta \ ;\
p\ ;\ \frac{\varepsilon }{2\gamma }).$

2$^{\grave{e}me}cas$ $1<q$ c'est-\`{a}-dire que $q^{\prime }<+\infty .$

Consid\`{e}rons un \'{e}l\'{e}ment $y$ de $G$. 
\begin{equation*}
\begin{array}{lll}
\left\| g\chi _{_{G\backslash B_{\left( e,\varepsilon \right) }}}\chi
_{_{yB_{\left( e,r\right) }}}\right\| _{q^{\prime }} & = & \left(
\int_{G}\left| g(x)\right| ^{q^{\prime }}\chi _{_{G\backslash B_{\left(
e,\varepsilon \right) }}}\left( x\right) \chi _{_{yB_{\left( e,r\right)
}}}\left( x\right) d\lambda (x)\right) ^{\frac{1}{q^{\prime }}} \\ 
& = & \left( \int_{G}\left| g(x)\right| ^{q^{\prime }}\chi _{_{G\backslash
B_{\left( e,\varepsilon \right) }}}\left( x\right) \chi _{_{B_{\left(
e,r\right) }}}\left( y^{-1}x\right) d\lambda (x)\right) ^{\frac{1}{q^{\prime
}}} \\ 
& = & \left( \int_{G}\left| g(yu)\right| ^{q^{\prime }}\chi _{_{G\backslash
B_{\left( e,\varepsilon \right) }}}\left( yu\right) \chi _{_{B_{\left(
e,r\right) }}}\left( u\right) d\lambda (u)\right) ^{\frac{1}{q^{\prime }}}.%
\end{array}
\end{equation*}
Nous avons pour tout \'{e}l\'{e}ment $u$ de $B_{\left( e,r\right) },$%
\begin{equation*}
yu\in G\backslash B_{\left( e,\varepsilon \right) }\Longrightarrow
\varepsilon <\left| yu\right| <\gamma \left( \left| y\right| +r\right)
\Longrightarrow \frac{\varepsilon }{\gamma }-r<\left| y\right|
\Longrightarrow \frac{\varepsilon }{2\gamma }<\left| y\right| . 
\end{equation*}
Donc, 
\begin{equation*}
\left\| g\chi _{_{G\backslash B_{\left( e,\varepsilon \right) }}}\chi
_{_{yB_{\left( e,r\right) }}}\right\| _{q^{\prime }}\leq \left(
\int_{G}\left| g(yu)\right| ^{q^{\prime }}\chi _{_{G\backslash B_{\left( e,%
\frac{\varepsilon }{2\gamma }\right) }}}\left( y\right) \chi _{_{B_{\left(
e,r\right) }}}\left( u\right) d\lambda (u)\right) ^{\frac{1}{q^{\prime }}}. 
\end{equation*}
Or, pour tout \'{e}l\'{e}ment $u$ de $B_{\left( e,r\right) },$ nous avons 
\begin{equation*}
\left| y\right| \leq \gamma \left( \left| yu\right| +\left| u\right| \right)
<\gamma \left| yu\right| +\gamma r<\gamma \left| yu\right| +\gamma 
\end{equation*}
\begin{equation*}
g(yu)=\left( \gamma \left| yu\right| +\gamma \right) ^{-\frac{\rho }{\theta
^{\prime }}}<\left| y\right| ^{-\frac{\rho }{\theta ^{\prime }}} 
\end{equation*}
Par suite 
\begin{equation*}
\left\| g\chi _{_{G\backslash B_{\left( e,\varepsilon \right) }}}\chi
_{_{yB_{\left( e,r\right) }}}\right\| _{q^{\prime }}\leq \lambda \left(
B_{\left( e,r\right) }\right) ^{\frac{1}{q^{\prime }}}\left| y\right| ^{-%
\frac{\rho }{\theta ^{\prime }}}\chi _{_{G\backslash B_{\left( e,\frac{%
\varepsilon }{2\gamma }\right) }}}\left( y\right) . 
\end{equation*}
Ainsi

$_{B_{\left( e,r\right) }}\left\| g\chi _{_{G\backslash B_{\left(
e,\varepsilon \right) }}}\right\| _{q^{\prime },p^{\prime }}\leq k(\theta \
;\ p\ ;\ \frac{\varepsilon }{2\gamma })\lambda \left( B_{\left( e,r\right)
}\right) ^{\frac{1}{q^{\prime }}}$
\end{proof}

\bigskip

\begin{proposition}
\label{prop: p<a}Supposons que $q,p$ et $\alpha $ sont des \'{e}l\'{e}ments
de $\left[ 1\ ;\ +\infty \right] $ tels que $p<\alpha .$

Alors $\left( L^{q},L^{p}\right) ^{\alpha }(G)=\left\{ O\right\} .$
\end{proposition}

\begin{proof}
Soit $f$ un \'{e}l\'{e}ment de $\left( L^{q},L^{p}\right) ^{\alpha }(G).$

Supposons que :

$\bullet \ \theta ,\varepsilon $ et $r$ sont des r\'{e}els tels que : $%
p<\theta <\alpha ,\ 0<\varepsilon <2\gamma $ et $0<r<\dfrac{\varepsilon }{%
2\gamma }$

$\bullet \ \forall x\in G\quad g(x)=\left( \gamma +\gamma \left| x\right|
\right) ^{-\frac{\rho }{\theta ^{\prime }}}$

Nous avons 
\begin{equation*}
\left\| fg\chi _{_{G\backslash B_{\left( e,\varepsilon \right) }}}\right\|
_{1}=\left\| fg\chi _{_{G\backslash B_{\left( e,\varepsilon \right)
}}}\right\| _{1,1}^{\pi _{r}}\leq \left\| f\right\| _{q,p}^{\pi _{r}}\left\|
g\chi _{_{G\backslash B_{\left( e,\varepsilon \right) }}}\right\|
_{q^{\prime },p^{\prime }}^{\pi _{r}}. 
\end{equation*}
Or d'apr\`{e}s la proposition \ref{PROP: équivalence BusbyetB( )qp}, il
existe un nombre r\'{e}el $C$ d\'{e}pendant uniquement de $G$ tel que 
\begin{equation*}
\left\| g\chi _{_{G\backslash B_{\left( e,\varepsilon \right) }}}\right\|
_{q^{\prime },p^{\prime }}^{\pi _{r}}\leq C\ _{B_{\left( e,r\right)
}}\left\| g\chi _{_{G\backslash B_{\left( e,\varepsilon \right) }}}\right\|
_{q^{\prime },p^{\prime }}\lambda \left( B_{\left( e,r\right) }\right) ^{-%
\frac{1}{p^{\prime }}}. 
\end{equation*}
En plus d'apr\`{e}s le lemme \ref{Lem: B()q'p'<kr}, 
\begin{equation*}
_{B_{\left( e,r\right) }}\left\| g\chi _{_{G\backslash B_{\left(
e,\varepsilon \right) }}}\right\| _{q^{\prime },p^{\prime }}\leq k(\theta \
;\ p\ ;\ \frac{\varepsilon }{2\gamma })\lambda \left( B_{\left( e,r\right)
}\right) ^{\frac{1}{q^{\prime }}} 
\end{equation*}
Par cons\'{e}quent, 
\begin{equation*}
\left\| fg\chi _{_{G\backslash B_{\left( e,\varepsilon \right) }}}\right\|
_{1}\leq Ck(\theta \ ;\ p\ ;\ \frac{\varepsilon }{2\gamma })\lambda \left(
B_{\left( e,r\right) }\right) ^{\frac{1}{p}-\frac{1}{\alpha }}\left\|
f\right\| _{q,p,\alpha }. 
\end{equation*}
$\frac{1}{p}-\frac{1}{\alpha }$ \'{e}tant strictement positif, nous
obtenons 
\begin{equation*}
\left\| fg\chi _{_{G\backslash B_{\left( e,\varepsilon \right) }}}\right\|
_{1}=0, 
\end{equation*}
en faisant tendre $r$ vers 0 dans l'in\'{e}galit\'{e} pr\'{e}c\'{e}dente.

Or pour tout \'{e}l\'{e}ment $x$ de $G\backslash B_{\left( e,\varepsilon
\right) },\quad g(x)\neq 0.$ Donc, 
\begin{equation*}
f\chi _{_{G\backslash B_{\left( e,\varepsilon \right) }}}=0. 
\end{equation*}
En faisant tendre $\varepsilon $ vers 0, nous obtenons \ $f=0. $
\end{proof}

\bigskip

Les propositions \ref{prop:cas a<q} et \ref{prop: p<a} montrent que l'espace 
$\left( L^{q},L^{p}\right) ^{\alpha }(G)$ est non trivial uniquement lorsque 
$q\leq \alpha \leq p$.

\begin{proposition}
\ Soient $\left( q_{1},p_{1},\alpha _{1}\right) $\ et $\left(
q_{2},p_{2},\alpha _{2}\right) $\ deux \'{e}l\'{e}ments de $\left[ 1\ ;\
+\infty \right] ^{3}$\ tels que $q_{1}\leq \alpha _{1}$\ et $q_{2}\leq
\alpha _{2}$. Si 

\begin{equation*}
\frac{1}{q_{1}}+\frac{1}{q_{2}}=\frac{1}{q}\leq 1,\ \frac{1}{p_{1}}+%
\frac{1}{p_{2}}=\frac{1}{p}\leq 1\ \text{et}\ \frac{1}{\alpha _{1}}+%
\frac{1}{\alpha _{2}}=\frac{1}{\alpha }, 
\end{equation*}

alors, pour $f$\ et $g$\ \'{e}l\'{e}ments de $L_{0}(G)$, nous avons :
 
\begin{equation*}
\left\| fg\right\| _{q,p,\alpha }\leq \left\| f\right\| _{q_{1},p_{1},\alpha
_{1}}\left\| g\right\| _{q_{2},p_{2},\alpha _{2}}. 
\end{equation*}
\end{proposition}

\begin{proof}
Elle est imm\'{e}diate \`{a} partir des d\'{e}finitions des normes $\left\|
{}\right\| _{q,p}^{\pi }$ et $\left\| {}\right\| _{q,p,\alpha }$.
\end{proof}

\section{LIEN AVEC LES ESPACES DE LEBESGUE ET DE LORENTZ.}

\begin{remark}
\label{REM:incl()qpa}Soient $p_{1},p_{2},q_{1},q_{2}$ et $\alpha $ des \'{e}l%
\'{e}ments de $\left[ 1\ ;\ +\infty \right] $ tels que $q_{1}\leq q_{2}\leq
\alpha \leq p_{2}\leq p_{1}.$ Alors 
\begin{equation*}
\left\| f\right\| _{q_{1},p_{1},\alpha }\leq \left( \dfrac{1}{2\gamma }%
\right) ^{\rho \left( \frac{1}{q_{1}}-\frac{1}{q_{2}}\right) }\left\|
f\right\| _{q_{2},p_{2},\alpha }\ 
\end{equation*}
et 
\begin{equation*}
\left\| f\right\| _{q_{1},q_{1},q_{1}}=\left\| f\right\| _{q_{1}}, 
\end{equation*}
pour tout \'{e}l\'{e}ment $f$ de $L_{0}(G).$
\end{remark}

\begin{proof}
Elle est imm\'{e}diate \`{a} partir de la d\'{e}finition de la norme $%
\left\| {}\right\| _{q,p,\alpha }.$
\end{proof}

\begin{proposition}
\label{Prop:limBr()q+inf=()q}Soit $q$\ un \'{e}l\'{e}ment de $\left[ 1\ ;\
+\infty \right[ $\ et $f$\ un \'{e}l\'{e}ment de \ $\left( L^{q},L^{+\infty
}\right) \left( G\right) $.\ 
\begin{equation*}
\underset{r\rightarrow +\infty }{\lim }\text{ }_{B_{\left( e,r\right)
}}\left\| f\right\| _{q,+\infty }=\left\| f\right\| _{q}=\underset{r>0}{\sup 
}\text{ }_{B_{\left( e,r\right) }}\left\| f\right\| _{q,+\infty } 
\end{equation*}
\end{proposition}

\begin{proof}
Soient $f$\ un \'{e}l\'{e}ment de $\left( L^{q},L^{+\infty }\right) \left(
G\right) $, $y$\ un \'{e}l\'{e}ment de $G$\ et $r$\ un r\'{e}el strictement
positif. 
\begin{equation*}
\left\| f\chi _{_{yB_{\left( e,r\right) }}}\right\| _{q}=\left[
\int_{yB_{\left( e,r\right) }}\left| f(t)\right| ^{q}d\lambda (t)\right] ^{%
\frac{1}{q}}\leq \left[ \int_{G}\left| f(t)\right| ^{q}d\lambda (t)\right] ^{%
\frac{1}{q}}=\left\| f\right\| _{q} 
\end{equation*}

\ Donc, pour tout r\'{e}el $\ r>0,$\ 
\begin{equation*}
\ _{B_{\left( e,r\right) }}\left\| f\right\| _{q,+\infty }=\underset{y\in G}{%
\sup }\left\| f\chi _{_{yB_{\left( e,r\right) }}}\right\| _{q}\leq \left\|
f\right\| _{q}. 
\end{equation*}

\ Par cons\'{e}quent, 
\begin{equation*}
\underset{r>0}{\sup }\text{ }_{B_{\left( e,r\right) }}\left\| f\right\|
_{q,+\infty }=\underset{r\rightarrow +\infty }{\lim }\text{ }_{B_{\left(
e,r\right) }}\left\| f\right\| _{q,+\infty }\leq \left\| f\right\| _{q}. 
\end{equation*}

Si $\underset{r>0}{\sup }$\ $_{B_{\left( e,r\right) }}\left\| f\right\|
_{q,+\infty }=+\infty $,\ alors l'\'{e}galit\'{e} s'ensuit.

Supposons $:$%
\begin{equation*}
\underset{r>0}{\sup }\text{ }_{B_{\left( e,r\right) }}\left\| f\right\|
_{q,+\infty }=M<+\infty . 
\end{equation*}

Pour tout r\'{e}el $r>0$, nous avons :\ 
\begin{equation*}
\ _{B_{\left( e,r\right) }}\left\| f\right\| _{q,+\infty }\leq M\text{.} 
\end{equation*}

\ Cela signifie que pour $\lambda -$presque tout $y$\ dans $G$\ et pour tout
r\'{e}el $r>0$,\ 
\begin{equation*}
\int_{B_{\left( e,r\right) }}\left| \text{ }_{y^{-1}}f(t)\right|
^{q}d\lambda (t)\leq M^{q}. 
\end{equation*}

D'o\`{u} 
\begin{equation*}
\left\| f\right\| _{q}\leq M=\underset{r>0}{\sup }\ _{B_{\left( e,r\right)
}}\left\| f\right\| _{q,+\infty }, 
\end{equation*}
puisque\ 

$\qquad \qquad \qquad \underset{r>0}{\sup }\int_{B_{\left( e,r\right)
}}\left| \text{ }_{y^{-1}}f(t)\right| ^{q}d\lambda (t)=\left\| \text{ }%
_{y^{-1}}f\right\| _{q}^{q}=\left\| f\right\| _{q}^{q}. $
\end{proof}

\bigskip

\begin{corollary}
\label{cor: La(G) inclu (Lq,Lp) a(G)}Pour tous \'{e}l\'{e}ments $p,q$\ et $%
\alpha $\ de $\left[ 1\ ;\ +\infty \right] $\ tels que $q\leq \alpha \leq p,$%
\begin{equation}
\left\| f\right\| _{q,p,\alpha }\leq \left\| f\right\| _{\alpha }\ \forall
f\in L_{0}(G).  \label{4}
\end{equation}
Par suite, $L^{\alpha }(G)\subset \left( L^{q},L^{p}\right) ^{\alpha }(G).$%
\end{corollary}

\bigskip

\begin{proposition}
\label{PROP : L qqin=Lq}Soit $q$\ un \'{e}l\'{e}ment de $\left[ 1\ ;\
+\infty \right] $. Il existe une constante $K=K\left( \gamma ,\rho \right) $
telle que pour tout \'{e}l\'{e}ment $f$\ de $L_{0}(G),$ 
\begin{equation*}
K\left\| f\right\| _{q}\leq \left\| f\right\| _{q,+\infty ,q}\leq \left\|
f\right\| _{q}. 
\end{equation*}
Par suite, $L^{q}(G)=\left( L^{q},L^{+\infty }\right) ^{q}(G).$
\end{proposition}

\begin{proof}
\ Soit $f$\ un \'{e}l\'{e}ment de $L_{0}(G)$. Nous avons $\left\| f\right\|
_{q,+\infty ,q}\leq \left\| f\right\| _{q}$, d'apr\`{e}s la relation $\left( %
\ref{4}\right) .$

Par ailleurs, pour tout r\'{e}el $r>0$, nous avons d'apr\`{e}s la
proposition \ref{PROP: équivalence BusbyetB( )qp},
 
\begin{equation*}
_{B_{\left( e,r\right) }}\left\| f\right\| _{q,+\infty }\leq C\left\|
f\right\| _{q,+\infty }^{\pi _{r}}\ ; 
\end{equation*}

ce qui signifie que
 
\begin{equation*}
\underset{r>0}{\sup }\text{ }_{B_{\left( e,r\right) }}\left\| f\right\|
_{q,+\infty }\leq C\underset{r>0}{\sup }\left\| f\right\| _{q,+\infty }^{\pi
_{r}}=C\left\| f\right\| _{q,+\infty ,q}. 
\end{equation*}

D'o\`{u} 

\begin{equation*}
K\left\| f\right\| _{q}\leq \left\| f\right\| _{q,+\infty ,q}\leq \left\|
f\right\| _{q}, 
\end{equation*}

puisque\ $\left\| f\right\| _{q}=\underset{r>0}{\sup }$\ $_{B_{\left(
e,r\right) }}\left\| f\right\| _{q,+\infty }. $
\end{proof}

\bigskip

\begin{proposition}
\label{PROP: Lqpq=Lq}Soient $p$\ et $q$\ deux \'{e}l\'{e}ments de $\left[ 1\
;\ +\infty \right] $, avec $1\leq q\leq p\leq +\infty $. Nous avons : 

\begin{equation*}
\left( L^{q},L^{p}\right) ^{q}(G)=L^{q}(G). 
\end{equation*}

Plus pr\'{e}cisement si $1\leq q\leq p<+\infty $, alors il existe une
constante $C=C\left( \gamma ,\rho ,p,q\right) $ 

\begin{equation*}
\ \left\| f\right\| _{q,p,q}\leq \left\| f\right\| _{q}\leq C\left\|
f\right\| _{q,p,q} 
\end{equation*}

\ pour tout \'{e}l\'{e}ment $f$\ de $L_{0}(G)$.
\end{proposition}

\begin{proof}
\ Soit $f$\ un \'{e}l\'{e}ment de $L_{0}(G)$.

\ Nous avons $\left\| f\right\| _{q,p,q}\leq \left\| f\right\| _{q},$ d'apr%
\`{e}s le corollaire \ref{cor: La(G) inclu (Lq,Lp) a(G)}$.$

\ Si $p=q$, alors il n'y a rien \`{a} d\'{e}montrer.

\ Si $q<p=+\infty $, alors nous retrouvons la proposition \ref{PROP : L
qqin=Lq}.

\ Supposons $q<p<+\infty $. Pour tout r\'{e}el $r>0$, 
\begin{equation*}
\begin{array}{ll}
\left\| f\chi _{_{B_{\left( e,r\right) }}}\right\| _{q}^{q} & 
=\int_{G}\left| f(x)\right| ^{q}\chi _{_{B_{\left( e,r\right) }}}(x)d\lambda
(x) \\ 
& =\underset{E\in \pi _{r}}{\sum }\int_{_{E}}\left| f(x)\right| ^{q}\chi
_{_{B_{\left( e,r\right) }}}(x)d\lambda (x) \\ 
& =\underset{E\in T_{e}}{\sum }\int_{G}\left| \left( f\chi _{_{E}}\right)
(x)\right| ^{q}\chi _{_{B_{\left( e,r\right) }}}(x)d\lambda (x) \\ 
& \leq \underset{E\in T_{e}}{\sum }\left\| f\chi _{_{E}}\right\|
_{q}^{q}\leq \left( 4\gamma ^{4}+3\gamma ^{2}\right) ^{\frac{\rho \left(
p-q\right) }{p}}\left[ \underset{E\in T_{e}}{\sum }\left\| f\chi
_{_{E}}\right\| _{q}^{p}\ \right] ^{\frac{q}{p}}.%
\end{array}
\end{equation*}

Donc, pour tout r\'{e}el strictement positif $r$, 
\begin{equation*}
\left\| f\chi _{_{B_{\left( e,r\right) }}}\right\| _{q}\leq \left( 4\gamma
^{4}+3\gamma ^{2}\right) ^{\frac{\rho \left( p-q\right) }{pq}}\left\|
f\right\| _{q,p}^{\pi _{r}}, 
\end{equation*}
ce qui signifie que\ 
\begin{equation*}
\ \left\| f\right\| _{q}\leq \left( 4\gamma ^{4}+3\gamma ^{2}\right) ^{\frac{%
\rho \left( p-q\right) }{pq}}\left\| f\right\| _{q,p,q}. 
\end{equation*}
$\ $

D'o\`{u},$\qquad \qquad \qquad \qquad \left\| f\right\| _{q,p,q}\leq \left\|
f\right\| _{q}\leq \left( 4\gamma ^{4}+3\gamma ^{2}\right) ^{\frac{\rho
\left( p-q\right) }{pq}}\left\| f\right\| _{q,p,q}.$
\end{proof}

\begin{proposition}
\label{PROP: Lqpp=Lp}Soient $p$\ et $q$\ deux \'{e}l\'{e}ments de $\left[ 1\
;\ +\infty \right] $\ tels que $q\leq p$ et $f$\ un \'{e}l\'{e}ment de $%
L_{0}(G)$. Alors il existe une constante $C=C\left( \gamma ,\rho ,p,q\right) 
$ telle que 
\begin{equation*}
\left\| f\right\| _{q,p,p}\leq \left\| f\right\| _{p}\leq C\left\| f\right\|
_{q,p,p}. 
\end{equation*}
Par suite, $L^{p}(G)=\left( L^{q},L^{p}\right) ^{p}(G).$
\end{proposition}

\begin{proof}
Soit $f$\ un \'{e}l\'{e}ment de $L_{0}(G)$

\ Nous avons $\left\| f\right\| _{q,p,p}\leq \left\| f\right\| _{p},$ d'apr%
\`{e}s le corollaire \ref{cor: La(G) inclu (Lq,Lp) a(G)}$.$\ 

Si $p=q$, ou $\left\| f\right\| _{q,p,p}=+\infty $,\ alors nous avons\ $%
\left\| f\right\| _{q,p,p}=\left\| f\right\| _{p}.$

Supposons $p\neq q$\ et $\left\| f\right\| _{q,p,p}<+\infty .$

Pour tout r\'{e}el $r>0,$\ nous avons :\ 
\begin{equation*}
\left\{ 
\begin{array}{l}
\lambda \left( B_{\left( e,r\right) }\right) ^{-\frac{1}{q}}\left\|
f\right\| _{q,+\infty }^{\pi _{r}}\leq \left\| f\right\| _{q,+\infty
,+\infty }<+\infty \text{ } \\ 
\ _{B_{\left( e,r\right) }}\left\| f\right\| _{q,+\infty }\leq \left(
4\gamma ^{5}+3\gamma ^{3}+2\gamma ^{2}\right) ^{\rho }\left\| f\right\|
_{q,+\infty }^{\pi _{r}}%
\end{array}
\right. \text{ si\ }q<p=+\infty , 
\end{equation*}

et 
\begin{equation*}
\ _{B_{\left( e,r\right) }}\left\| f\right\| _{q,p}\leq \lambda \left(
B_{\left( e,r\right) }\right) ^{\frac{1}{q}}\left( 4\gamma ^{4}+3\gamma
^{2}\right) ^{\frac{\rho }{q}}\left( 4\gamma ^{5}+3\gamma ^{3}+2\gamma
^{2}\right) ^{\frac{\rho }{p}}\left\| f\right\| _{q,p,p}\ \text{si}\
q<p<+\infty . 
\end{equation*}

Posons :\ 
\begin{equation*}
f_{r}(x)=\left[ \dfrac{1}{\lambda \left( B_{\left( e,r\right) }\right) }%
\int_{B_{\left( x,r\right) }}\left| f(t)\right| ^{q}d\lambda (t)\right] ^{%
\frac{1}{q}}. 
\end{equation*}

Nous avons :\ 
\begin{equation*}
\underset{r\rightarrow 0}{\lim }f_{r}(x)=\left| f(x)\right| \leq \left(
4\gamma ^{5}+3\gamma ^{3}+2\gamma ^{2}\right) ^{\rho }\left\| f\right\|
_{q,+\infty ,+\infty }, 
\end{equation*}
pour $\lambda -$presque tout $x$ dans $G$. Par cons\'{e}quent, 
\begin{equation*}
\left\| f\right\| _{+\infty }\leq \left( 4\gamma ^{5}+3\gamma ^{3}+2\gamma
^{2}\right) ^{\rho }\left\| f\right\| _{q,+\infty ,+\infty }. 
\end{equation*}

Par ailleurs, 
\begin{equation*}
\left[ \int_{G}f_{r}^{\text{ }p}(x)d\lambda (x)\right] ^{\frac{1}{p}}\leq
\left( 4\gamma ^{4}+3\gamma ^{2}\right) ^{\frac{\rho }{q}}\left( 4\gamma
^{5}+3\gamma ^{3}+2\gamma ^{2}\right) ^{\frac{\rho }{p}}\left\| f\right\|
_{q,p,p}. 
\end{equation*}
\ 

D'o\`{u}, d'apr\`{e}s le lemme de Fatou, $\left| f\right| ^{p}$est int\'{e}%
grable et 
\begin{equation*}
\left\| f\right\| _{p}\leq \left( 4\gamma ^{4}+3\gamma ^{2}\right) ^{\frac{%
\rho }{q}}\left( 4\gamma ^{5}+3\gamma ^{3}+2\gamma ^{2}\right) ^{\frac{\rho 
}{p}}\left\| f\right\| _{q,p,p}. 
\end{equation*}

Par suite $\left( L^{q},L^{p}\right) ^{p}(G)=L^{p}(G).$\bigskip
\end{proof}

\begin{proposition}
\label{PROP: ( )qpa< ( )*a} Soit $\left( q,p,\alpha \right) $\ un \'{e}l\'{e}%
ment de $\left[ 1\ ;\ +\infty \right] ^{3}$\ tels que $q<\alpha <p.$ Il
existe une constante $C=C(q,p,\alpha ,\rho )$\ telle que 

\begin{equation*}
\left\| f\right\| _{q,p,\alpha }\leq C\left\| f\right\| _{\alpha ,+\infty
}^{\ast },\ \forall f\in L_{0}(G). 
\end{equation*}

Par suite, $L^{\alpha ,+\infty }(G)\subset \left( L^{q},L^{p}\right)
^{\alpha }(G).$
\end{proposition}

\begin{proof}
Soit $f$\ un \'{e}l\'{e}ment de $L^{\alpha ,+\infty }(G)$.

\textbf{1}$^{er}$\textbf{\ cas }: supposons $p=+\infty $.

Pour tout r\'{e}el $r>0$\ et pour $\lambda -$presque tout $x$\ dans $G$,
nous avons :\ 
\begin{equation*}
\left( \left| f\right| ^{q}\ast \chi _{_{B_{\left( e,r\right) }}}\right)
(x)=\int_{G}\left| f(t)\right| ^{q}\chi _{_{B_{\left( x,r\right)
}}}(t)d\lambda (t)=\left\| f\chi _{_{B_{\left( x,r\right) }}}\right\|
_{q}^{q}. 
\end{equation*}

Puisque $1\leq q<\alpha $, nous avons d'apr\`{e}s la condition dite de
Kolmogorov \cite{3},\ 
\begin{equation*}
\left\| f\chi _{_{B_{\left( x,r\right) }}}\right\| _{q}\leq \left( \frac{%
\alpha }{\alpha -q}\right) ^{\frac{1}{q}}\left\| f\right\| _{\alpha ,+\infty
}^{\ast }\lambda \left( B_{\left( e,r\right) }\right) ^{\frac{1}{q}-\frac{1}{%
\alpha }}. 
\end{equation*}

\ Par cons\'{e}quent, pour $\lambda -$presque tout $x$\ dans $G,$\ 
\begin{equation*}
\left( \left| f\right| ^{q}\ast \chi _{_{B_{\left( e,r\right) }}}\right)
(x)\leq \left[ \left( \frac{\alpha }{\alpha -q}\right) ^{\frac{1}{q}}\left\|
f\right\| _{\alpha ,+\infty }^{\ast }\lambda \left( B_{\left( e,r\right)
}\right) ^{\frac{1}{q}-\frac{1}{\alpha }}\right] ^{q}. 
\end{equation*}

Donc, 
\begin{equation*}
\left\| \left( \left| f\right| ^{q}\ast \chi _{_{B_{\left( e,r\right)
}}}\right) ^{\frac{1}{q}}\right\| _{+\infty }\leq \left( \frac{\alpha }{%
\alpha -q}\right) ^{\frac{1}{q}}\left\| f\right\| _{\alpha ,+\infty }^{\ast
}\lambda \left( B_{\left( e,r\right) }\right) ^{\frac{1}{q}-\frac{1}{\alpha }%
}\ ; 
\end{equation*}
Ce qui peut encore s'\'{e}crire :\ 
\begin{equation*}
\lambda \left( B_{\left( e,r\right) }\right) ^{\frac{1}{\alpha }-\frac{1}{q}}%
\text{ }_{B_{\left( e,r\right) }}\left\| f\right\| _{q,+\infty }\leq \left( 
\dfrac{\alpha }{\alpha -q}\right) ^{\frac{1}{q}}\left\| f\right\| _{\alpha
,+\infty }^{\ast }. 
\end{equation*}

D'o\`{u},\ 
\begin{equation*}
\left\| f\right\| _{q,+\infty ,\alpha }\leq \left( \dfrac{\alpha }{\alpha -q}%
\right) ^{\frac{1}{q}}\left\| f\right\| _{\alpha ,+\infty }^{\ast }. 
\end{equation*}

\textbf{2}$^{eme}$\textbf{\ cas }: supposons que $p<+\infty $.

Posons $\qquad \beta =\left( 1+\dfrac{q}{p}-\dfrac{q}{\alpha }\right) ^{-1}$.

Alors $1<\beta <+\infty $, $1<\dfrac{\alpha }{q}<+\infty $\ et $\dfrac{q}{p}=%
\dfrac{1}{\beta }+\dfrac{q}{\alpha }-1$.

Puisque $\left| f\right| ^{q}$\ est un \'{e}l\'{e}ment de $L^{\frac{\alpha }{%
q},+\infty }(G)$\ et $\chi _{_{B_{\left( e,r\right) }}}$\ un \'{e}l\'{e}ment
de $L^{\beta }(G)$, 
\begin{equation*}
\left| f\right| ^{q}\ast \chi _{_{B_{\left( e,r\right) }}}\in L^{\frac{p}{q}%
}(G), 
\end{equation*}
et 
\begin{equation*}
\left\| \left| f\right| ^{q}\ast \chi _{_{B_{\left( e,r\right) }}}\right\| _{%
\frac{p}{q}}\leq C\left( \left\| f\right\| _{\alpha ,+\infty }^{\ast
}\right) ^{q}\lambda \left( B_{\left( e,r\right) }\right) ^{\frac{1}{\beta }%
}. 
\end{equation*}

Remarquons que$\qquad \dfrac{1}{\beta q}=\dfrac{1}{q}+\dfrac{1}{p}-\dfrac{1}{%
\alpha }$.

Nous avons, pour tout r\'{e}el $r>0$, 
\begin{equation*}
\left\| \left| f\right| ^{q}\ast \chi _{_{B_{\left( e,r\right) }}}\right\| _{%
\frac{p}{q}}^{\frac{1}{q}}\leq C\left\| f\right\| _{\alpha ,+\infty }^{\ast
}\lambda \left( B_{\left( e,r\right) }\right) ^{\frac{1}{q}+\frac{1}{p}-%
\frac{1}{\alpha }}. 
\end{equation*}

Ce qui peut encore s'\'{e}crire, pour tout r\'{e}el $r>0$, 
\begin{equation*}
\lambda \left( B_{\left( e,r\right) }\right) ^{\frac{1}{\alpha }-\frac{1}{p}-%
\frac{1}{q}}\text{ }_{B_{\left( e,r\right) }}\left\| f\right\| _{q,p}\leq
C\left\| f\right\| _{\alpha ,+\infty }^{\ast }. 
\end{equation*}

D'o\`{u},

$\qquad \qquad \qquad \qquad \left\| f\right\| _{q,p,\alpha }\leq C^{\prime }%
\underset{r>0}{\sup }\lambda \left( B_{\left( e,r\right) }\right) ^{\frac{1}{%
\alpha }-\frac{1}{p}-\frac{1}{q}}$ $_{B_{\left( e,r\right) }}\left\|
f\right\| _{q,p}\leq C^{\prime }C\left\| f\right\| _{\alpha ,+\infty }^{\ast
}$
\end{proof}

En fait, l'inclusion de la proposition pr\'{e}c\'{e}dente est stricte.

\begin{proposition}
\label{PROP : element de Lqpa} Soit $\left( q,p,\alpha \right) $\ un \'{e}l%
\'{e}ment de $\left[ 1\ ;\ +\infty \right] ^{3}$\ tel que $1\leq q<\alpha
<p. $ Consid\'{e}rons dans $G$ une famille de boules 
\begin{equation*}
\left\{ B_{\left( x_{j}^{n},2^{-n-1}\right) }\text{ }/\text{ }1\leq j\leq
E\left( 2^{\rho \left( n+1\right) }\right) +1\ ;\ n\in \mathbb{N}^{\ast
}\right\} , 
\end{equation*}
o\`{u} $E\left( 2^{\rho \left( n+1\right) }\right) $ d\'{e}signe la partie
enti\`{e}re de $2^{\rho \left( n+1\right) },$ v\'{e}rifiant les conditions
suivantes\ :

i) $\left| \left( x_{k}^{1}\right) ^{-1}x_{j}^{1}\right| >\gamma \left(
1+\gamma +2\gamma ^{2}\right) 2^{\frac{2q}{\alpha -q}},\forall \left(
k,j\right) \in \left\{ 1;2;3;\ldots ;E\left( 2^{2\rho }\right) +1\right\}
^{2},$ \newline
$j\neq k\ ;$

ii) pour tout entier $n>1,$

\qquad $\vartriangleright \left| \left( x_{k}^{n}\right)
^{-1}x_{j}^{n}\right| >\gamma \left( 1+\gamma +2\gamma ^{2}\right) 2^{\frac{%
\left( n+1\right) q}{\alpha -q}},\ \forall \left( k,j\right) \in \left\{
1;2;3;\ldots ;E\left( 2^{\rho \left( n+1\right) }\right) +1\right\} ^{2},$ 
\newline
$j\neq k,$ \newline
et

\qquad $\vartriangleright \left| \left( x_{k}^{n}\right)
^{-1}x_{j}^{m}\right| >2\gamma ^{2}\left| 2^{2^{\rho \left( n+1\right)
}+1}+2^{k-n}+2^{-n-1}-2^{2^{\rho \left( m+1\right)
}+1}-2^{j-m}-2^{-m-1}\right| ,$

$\forall \left( m,k,j\right) \in \left\{ 1;2;3;\ldots ;n\right\} \times
\left\{ 1;2;3;\ldots ;E\left( 2^{\rho \left( n+1\right) }\right) +1\right\}
\times \left\{ 1;2;3;\ldots ;E\left( 2^{\rho \left( m+1\right) }\right)
+1\right\} ,$\newline
avec $m\neq n$ ou $j\neq k.$

Posons : 
\begin{equation*}
E_{n}=\overset{E\left( 2^{\rho \left( n+1\right) }\right) +1}{\underset{j=1}{%
\cup }}B_{\left( x_{j}^{n},2^{-n-1}\right) }, 
\end{equation*}
\begin{equation*}
E=\underset{n\geq 1}{\cup }E_{n}, 
\end{equation*}
et\qquad\ $f=\chi _{_{E}}.$

Alors, $f$ est un \'{e}l\'{e}ment de $\left( L^{q},L^{p}\right) ^{\alpha
}(G) $ qui n'est pas dans $L^{\alpha ,+\infty }(G).$
\end{proposition}

\begin{proof}
La famille $\left\{ B_{\left( x_{j}^{n},2^{-n-1}\right) }\text{ }/\text{ }%
1\leq j\leq E\left( 2^{\rho \left( n+1\right) }\right) +1\ ;\text{ }n\in 
\mathbb{N}^{\ast }\right\} ,$ est compos\'{e}e de boules deux \`{a} deux
disjointes.

En effet, soient $m$ et $n$ deux entiers naturels non nuls, tels que $m\geq
n,$

$\left( k,j\right) \in \left\{ 1;2;3;\ldots ;E\left( 2^{\rho \left(
n+1\right) }\right) +1\right\} \times \left\{ 1;2;3;\ldots ;E\left( 2^{\rho
\left( m+1\right) }\right) +1\right\} $ .

Si $\ a$ est un \'{e}l\'{e}ment de $B_{\left( x_{j}^{n},2^{-n-1}\right)
}\cap B_{\left( x_{k}^{m},2^{-m-1}\right) }$ avec $j\neq k,$ alors nous
avons\ : 
\begin{equation*}
\begin{array}{ll}
\gamma 2^{-n} & <2\gamma ^{2}2^{-n}<2\gamma ^{2}\left( 2^{2^{\rho \left(
m+1\right) }+1}+2^{j-m}+2^{-m-1}-2^{2^{\rho \left( n+1\right)
}+1}-2^{k-n}-2^{-n-1}\right) \\ 
& <\left| \left( x_{k}^{n}\right) ^{-1}x_{j}^{m}\right| \leq \gamma \left(
\left| \left( x_{k}^{n}\right) ^{-1}a\right| +\left| a^{-1}x_{j}^{m}\right|
\right) <\gamma 2^{-n}\ ;%
\end{array}
\end{equation*}
ce qui est impossible.

Par cons\'{e}quent $\left( E_{n}\right) _{n\geq 1}$ est une famille
d'ensembles deux \`{a} deux disjoints.

Ainsi, nous avons 
\begin{equation*}
\lambda \left( E\right) =\underset{n\geq 1}{\sum }\lambda \left(
E_{n}\right) =\underset{n\geq 1}{\sum }\underset{j=1}{\overset{E\left(
2^{\rho \left( n+1\right) }\right) +1}{\sum }}\lambda \left( B_{\left(
x_{j}^{n},2^{-n-1}\right) }\right) =\underset{n\geq 1}{\sum }\left[ E\left(
2^{\rho \left( n+1\right) }\right) +1\right] \left( 2^{-n-1}\right) ^{\rho
}=+\infty . 
\end{equation*}

Par cons\'{e}quent, 
\begin{equation*}
f^{\ast }(t)=\inf \left\{ s>0,\lambda \left( \left\{ x\in G\text{ }/\text{ }%
\left| f(x)\right| >s\right\} \right) \leq t\right\} =1,\forall t\in \mathbb{%
R}_{+}^{\ast }. 
\end{equation*}

D'o\`{u}, 
\begin{equation*}
\left\| f\right\| _{\alpha ,+\infty }^{\ast }=\underset{t>0}{\sup }t^{\frac{1%
}{\alpha }}f^{\ast }(t)=\underset{t>0}{\sup }t^{\frac{1}{\alpha }}=+\infty . 
\end{equation*}

Par ailleurs, pour tout r\'{e}el $r>0,$ nous avons\ : 
\begin{equation*}
_{B_{\left( e,r\right) }}\left\| \chi _{_{E_{n}}}\right\| _{q,p}\leq \left[ 
\underset{j=1}{\overset{E\left( 2^{\rho \left( n+1\right) }\right) +1}{\sum }%
}\int_{B_{\left( x_{j}^{n},\gamma 2^{-n-1}(1+r2^{n+1})\right) }}\left(
\lambda \left( E_{n}\cap B_{\left( x,r\right) }\right) \right) ^{\frac{p}{q}%
}d\lambda (x)\right] ^{\frac{1}{p}},\forall n\in \mathbb{N}^{\ast }\ ; 
\end{equation*}
car, si $B_{\left( x_{j}^{n},2^{-n-1}\right) }\cap B_{\left( x,r\right)
}\neq \emptyset $ alors $\left| \left( x_{j}^{n}\right) ^{-1}x\right|
<\gamma \left( 2^{-n-1}+r\right) =\gamma 2^{-n-1}\left( 1+r2^{n+1}\right) .$

Soit $n$ un \'{e}l\'{e}ment de $\mathbb{N}^{\ast }.$

\qquad \textbf{a)} Supposons $r=2^{-n-1-\ell },$ avec $\ell $ un entier
naturel.

Alors, 
\begin{equation*}
\begin{array}{l}
\lambda \left( B_{\left( e,r\right) }\right) ^{\frac{1}{\alpha }-\frac{1}{q}-%
\frac{1}{p}}\ _{B_{\left( e,r\right) }}\left\| \chi _{_{E_{n}}}\right\|
_{q,p} \\ 
\leq r^{\rho \left( \frac{1}{\alpha }-\frac{1}{q}-\frac{1}{p}\right) }\left[ 
\underset{j=1}{\overset{E\left( 2^{\rho \left( n+1\right) }\right) +1}{\sum }%
}\int_{B_{\left( x_{j}^{n},\gamma 2^{-n-1}(1+2^{-\ell }\right) }}\lambda
\left( B_{\left( x_{j}^{n},2^{-n-1}\right) }\cap B_{\left( x,r\right)
}\right) ^{\frac{p}{q}}d\lambda (x)\right] ^{\frac{1}{p}},%
\end{array}
\end{equation*}
car $B_{\left( x_{j}^{n},\gamma 2^{-n-1}(1+2^{-\ell })\right) }\cap
B_{\left( x_{k}^{n},\gamma 2^{-n-1}(1+2^{-\ell })\right) }=\emptyset $, si $%
1\leq j<k\leq E\left( 2^{\rho \left( n+1\right) }\right) +1.$

Donc, 
\begin{equation*}
\begin{array}{l}
\lambda \left( B_{\left( e,r\right) }\right) ^{\frac{1}{\alpha }-\frac{1}{q}-%
\frac{1}{p}}\ _{B_{\left( e,r\right) }}\left\| \chi _{_{E_{n}}}\right\|
_{q,p} \\ 
\begin{array}{l}
\leq \left( 2^{-n-1-\ell }\right) ^{\rho \left( \frac{1}{\alpha }-\frac{1}{q}%
-\frac{1}{p}\right) }\left[ \left( 2^{\rho \left( n+1\right) }+1\right)
\left( \gamma 2^{-n-1}(1+2^{-\ell })\right) ^{\rho }\left( 2^{-n-1-\ell
}\right) ^{\rho \frac{p}{q}}\right] ^{\frac{1}{p}} \\ 
\leq 2^{\left( \frac{1+\rho }{p}+\frac{\rho }{p}-\frac{\rho }{\alpha }%
\right) }\gamma ^{\frac{\rho }{p}}\left( 2^{-\rho \left( \frac{1}{\alpha }-%
\frac{1}{p}\right) }\right) ^{n}.%
\end{array}%
\end{array}
\end{equation*}
C'est \`{a} dire que 
\begin{equation*}
\lambda \left( B_{\left( e,r\right) }\right) ^{\frac{1}{\alpha }-1-\frac{1}{p%
}}\text{ }_{B_{\left( e,r\right) }}\left\| \chi _{_{E_{n}}}\right\|
_{1,p}\leq C_{1}\left( 2^{-\rho \left( \frac{1}{\alpha }-\frac{1}{p}\right)
}\right) ^{n}, 
\end{equation*}
avec $C_{1}=2^{\left( \frac{1+\rho }{p}+\frac{\rho }{p}-\frac{\rho }{\alpha }%
\right) }\gamma ^{\frac{\rho }{p}}.$

\qquad \textbf{b)} Supposons $r=2^{-n-1+\ell },$ avec $\ell $ un entier
naturel non nul.

\qquad \qquad Si $1\leq \ell \leq \frac{\left( n+1\right) \alpha }{\alpha -q%
},$ alors pour tout $x$ dans $G$, la boule $B_{\left( x,2^{-n-1+\ell
}\right) }$ ne peut rencontrer plus d'une boule $B_{\left(
x_{k}^{n},2^{-n-1}\right) }\ ;\ $car, 

\begin{equation*}
\left\{ 
\begin{array}{l}
y_{k}\in B_{\left( x,2^{-n-1+\ell }\right) }\cap B_{\left(
x_{k}^{n},2^{-n-1}\right) } \\ 
y_{k^{\prime }}\in B_{\left( x,2^{-n-1+\ell }\right) }\cap B_{\left(
x_{k^{\prime }}^{n},2^{-n-1}\right) }%
\end{array}
\right. \Rightarrow \left| \left( x_{k}^{n}\right) ^{-1}x_{k^{\prime
}}^{n}\right| <\gamma \left( 2\gamma ^{2}+\gamma +1\right) 2^{\frac{\left(
n+1\right) q}{\alpha -q}}\ ; 
\end{equation*}

Ce qui est impossible si $k\neq k^{\prime },$ au vu de l'hypoth\`{e}se.

Nous avons par cons\'{e}quent\ : 
\begin{equation*}
\begin{array}{l}
\lambda \left( B_{\left( e,r\right) }\right) ^{\frac{1}{\alpha }-\frac{1}{q}-%
\frac{1}{p}}\ _{B_{\left( e,r\right) }}\left\| \chi _{_{E_{n}}}\right\|
_{q,p} \\ 
\begin{array}{l}
\leq r^{\rho \left( \frac{1}{\alpha }-\frac{1}{q}-\frac{1}{p}\right) }\left[ 
\underset{j=1}{\overset{E\left( 2^{\rho \left( n+1\right) }\right) +1}{\sum }%
}\int_{B_{\left( x_{j}^{n},\gamma 2^{-n-1}(1+2^{\ell })\right) }}\left(
\lambda \left( E_{n}\cap B_{\left( x,r\right) }\right) \right) ^{\frac{p}{q}%
}d\lambda (x)\right] ^{\frac{1}{p}} \\ 
\leq \left( 2^{-n-1+\ell }\right) ^{\rho \left( \frac{1}{\alpha }-\frac{1}{q}%
-\frac{1}{p}\right) }\left[ \left( 2^{\rho \left( n+1\right) }+1\right)
\left( \gamma 2^{-n-1}(1+2^{\ell })\right) ^{\rho }\left( 2^{-n-1}\right)
^{\rho \frac{p}{q}}\right] ^{\frac{1}{p}} \\ 
\leq 2^{\rho \left( n+1\right) \left( \frac{1}{p}-\frac{1}{\alpha }\right)
}2^{\frac{\rho +1}{p}}\gamma ^{\frac{\rho }{p}}.%
\end{array}%
\end{array}
\end{equation*}

Ainsi, 
\begin{equation*}
\lambda \left( B_{\left( e,r\right) }\right) ^{\frac{1}{\alpha }-\frac{1}{q}-%
\frac{1}{p}}\text{ }_{B_{\left( e,r\right) }}\left\| \chi
_{_{E_{n}}}\right\| _{q,p}\leq \left( 2^{-\rho \left( \frac{1}{\alpha }-%
\frac{1}{p}\right) }\right) ^{n}2^{\frac{\rho +1}{p}+\frac{\rho }{p}-\frac{%
\rho }{\alpha }}\gamma ^{\frac{\rho }{p}}\text{ si \ }1\leq \ell \leq \frac{%
\left( n+1\right) \alpha }{\alpha -q}. 
\end{equation*}

\qquad \qquad Si $\ell >\dfrac{\left( n+1\right) \alpha }{\alpha -q},$ alors
nous avons\ : 
\begin{equation*}
\begin{array}{l}
\lambda \left( B_{\left( e,r\right) }\right) ^{\frac{1}{\alpha }-\frac{1}{q}-%
\frac{1}{p}}\ _{B_{\left( e,r\right) }}\left\| \chi _{_{E_{n}}}\right\|
_{q,p} \\ 
\begin{array}{l}
\leq r^{\rho \left( \frac{1}{\alpha }-\frac{1}{q}-\frac{1}{p}\right) }\left[ 
\underset{j=1}{\overset{E\left( 2^{\rho \left( n+1\right) }\right) +1}{\sum }%
}\int_{B_{\left( x_{j}^{n},\gamma 2^{-n-1}(1+2^{\ell })\right) }}\left(
\lambda \left( E_{n}\cap B_{\left( x,r\right) }\right) \right) ^{\frac{p}{q}%
}d\lambda (x)\right] ^{\frac{1}{p}} \\ 
\leq r^{\rho \left( \frac{1}{\alpha }-\frac{1}{q}-\frac{1}{p}\right) }\left[ 
\underset{j=1}{\overset{E\left( 2^{\rho \left( n+1\right) }\right) +1}{\sum }%
}\int_{B_{\left( x_{j}^{n},\gamma 2^{-n-1}(1+2^{\ell })\right) }}\lambda
\left( E_{n}\right) ^{\frac{p}{q}}d\lambda (x)\right] ^{\frac{1}{p}} \\ 
\leq \left( 2^{-n-1+\ell }\right) ^{\rho \left( \frac{1}{\alpha }-\frac{1}{q}%
-\frac{1}{p}\right) }\left[ \left( 2^{\rho \left( n+1\right) }+1\right)
\left( \gamma 2^{-n-1}(1+2^{\ell })\right) ^{\rho }\left[ \left( 2^{\rho
\left( n+1\right) }+1\right) \left( 2^{-n-1}\right) ^{\rho }\right] ^{\frac{p%
}{q}}\right] ^{\frac{1}{p}} \\ 
\leq 2^{-\rho \left( n+1\right) \left( \frac{1}{\alpha }-\frac{1}{p}\right)
\gamma ^{\frac{\rho }{p}}}2^{\frac{\rho +1}{p}+\frac{1}{q}}.%
\end{array}%
\end{array}
\end{equation*}

Ainsi, 
\begin{equation*}
\lambda \left( B_{\left( e,r\right) }\right) ^{\frac{1}{\alpha }-\frac{1}{q}-%
\frac{1}{p}}\text{ }_{B_{\left( e,r\right) }}\left\| \chi
_{_{E_{n}}}\right\| _{q,p}\ \leq \left( 2^{-\rho \left( \frac{1}{\alpha }-%
\frac{1}{p}\right) }\right) ^{n}\gamma ^{\frac{\rho }{p}}2^{\frac{\rho +1}{p}%
+\frac{1}{q}-\rho \left( \frac{1}{\alpha }-\frac{1}{p}\right) }\ \text{si }%
\ell >\frac{n+1}{1-\frac{1}{\alpha }}. 
\end{equation*}

Donc, pour $r=2^{-n-1+\ell },$ avec $\ell $ entier naturel non nul, nous
avons\ : 
\begin{equation*}
\lambda \left( B_{\left( e,r\right) }\right) ^{\frac{1}{\alpha }-\frac{1}{q}-%
\frac{1}{p}}\text{ }_{B_{\left( e,r\right) }}\left\| \chi
_{_{E_{n}}}\right\| _{q,p}\leq C_{2}\left( 2^{-\rho \left( \frac{1}{\alpha }-%
\frac{1}{p}\right) }\right) ^{n}\ ; 
\end{equation*}

avec $C_{2}=\max (\gamma ^{\frac{\rho }{p}}2^{\frac{\rho +1}{p}+\frac{1}{q}%
-\rho \left( \frac{1}{\alpha }-\frac{1}{p}\right) }\ ;\ 2^{\frac{\rho +1}{p}+%
\frac{\rho }{p}-\frac{\rho }{\alpha }}\gamma ^{\frac{\rho }{p}}).$

Nous pouvons donc dire \`{a} partir de a) et b) que si $r=2^{-m},\ m$ \'{e}%
tant un entier relatif, alors 
\begin{equation*}
\lambda \left( B_{\left( e,r\right) }\right) ^{\frac{1}{\alpha }-\frac{1}{q}-%
\frac{1}{p}}\text{ }_{B_{\left( e,r\right) }}\left\| \chi
_{_{E_{n}}}\right\| _{q,p}\leq C_{3}\left( 2^{-\rho \left( \frac{1}{\alpha }-%
\frac{1}{p}\right) }\right) ^{n}, 
\end{equation*}
avec $C_{3}=\max (C_{1},C_{2})$.

\textbf{c)} Supposons $r$ quelconque.

\qquad Il existe un unique entier relatif $m$ tel que $2^{-m-1}\leq r\leq
2^{-m}.$ par cons\'{e}quent, 
\begin{equation*}
\begin{array}{ll}
\lambda \left( B_{\left( e,r\right) }\right) ^{\frac{1}{\alpha }-\frac{1}{q}-%
\frac{1}{p}}\ _{B_{\left( e,r\right) }}\left\| \chi _{_{E_{n}}}\right\|
_{q,p} & \leq \lambda \left( B_{\left( e,2^{-m-1}\right) }\right) ^{\frac{1}{%
\alpha }-\frac{1}{q}-\frac{1}{p}}\ _{B_{\left( e,2^{-m}\right) }}\left\|
\chi _{_{E_{n}}}\right\| _{q,p} \\ 
& \leq 2^{-\rho \left( \frac{1}{\alpha }-\frac{1}{q}-\frac{1}{p}\right)
}C_{3}\left( 2^{-\rho \left( \frac{1}{\alpha }-\frac{1}{p}\right) }\right)
^{n}.%
\end{array}
\end{equation*}

Donc, pour tout r\'{e}el $r>0,$%
\begin{equation*}
\lambda \left( B_{\left( e,r\right) }\right) ^{\frac{1}{\alpha }-\frac{1}{q}-%
\frac{1}{p}}\text{ }_{B_{\left( e,r\right) }}\left\| \chi
_{_{E_{n}}}\right\| _{q,p}\leq C_{4}\left( 2^{-\rho \left( \frac{1}{\alpha }-%
\frac{1}{p}\right) }\right) ^{n}, 
\end{equation*}

avec\qquad\ $C_{4}=2^{-\rho \left( \frac{1}{\alpha }-\frac{1}{q}-\frac{1}{p}%
\right) }C_{3}.$

Ceci \'{e}tant vrai pour un entier quelconque $n$ de $\mathbb{N}^{\ast },$
nous avons : 
\begin{equation*}
\lambda \left( B_{\left( e,r\right) }\right) ^{\frac{1}{\alpha }-\frac{1}{q}-%
\frac{1}{p}}\ _{B_{\left( e,r\right) }}\left\| f\right\| _{q,p}\leq \underset%
{n\geq 1}{\sum }C_{4}\left( 2^{-\rho \left( \frac{1}{\alpha }-\frac{1}{p}%
\right) }\right) ^{n}=C_{4}\frac{2^{-\rho \left( \frac{1}{\alpha }-\frac{1}{%
p}\right) }}{1-2^{-\rho \left( \frac{1}{\alpha }-\frac{1}{p}\right) }},\
\forall r\in \mathbb{R}_{+}^{\ast }. 
\end{equation*}

D'o\`{u},\qquad $\left\| f\right\| _{q,p,\alpha }\leq \underset{r>0}{\sup }%
\lambda \left( B_{\left( e,r\right) }\right) ^{\frac{1}{\alpha }-\frac{1}{q}-%
\frac{1}{p}}$ $_{B_{\left( e,r\right) }}\left\| f\right\| _{q,p}\leq C_{4}%
\frac{2^{-\rho \left( \frac{1}{\alpha }-\frac{1}{p}\right) }}{1-2^{-\rho
\left( \frac{1}{\alpha }-\frac{1}{p}\right) }}<+\infty .$
\end{proof}

\bigskip

Les remarques faites par le ref\'{e}r\'{e} nous ont permis de pr\'{e}ciser
certains de nos r\'{e}sultats. Qu'il trouve ici l'expression de notre
profonde gratitude.


\bigskip

\begin{thebibliography}{1}
\bibitem[1]{1} Robert C. Busby and Harvey A. Smith, \textit{%
Product-convolution Operators and Mixed-Norm spaces}, Trans. AMS, vol 263,
Number 2, February 1981. 309-341.

\bibitem[2]{2} Ibrahim Fofana, \textit{Etude d'une classe de fonctions
contenant l'espace de Lorentz}, Afrika Mathematika (2) 1 (1988) 29-50.

\bibitem[3]{4} G. B. Folland, \textit{Subelliptic estimates and function
spaces on nilpotent Lie proups, }Arkiv mat. vol 13 (1975) n$%
{{}^\circ}%
2.$

\bibitem[4]{3} J.Garcia-Cuerva and J.L\ Rubio de Francia, \textit{Weighted
norm inequalities and related topics}, in ''North-Holland Math. Studies,''
vol. 116, North-Holland Amsterdam.(485)$\ $
\end{thebibliography}
\end{document}